\documentclass[journal]{IEEEtran}
\IEEEoverridecommandlockouts                              % This command is only
\usepackage{amsfonts}
\usepackage{amsmath}
\usepackage{amssymb}
\newtheorem{theorem}{\textbf{Theorem}}
\newtheorem{corollary}{\textbf{Corollary}}
\newtheorem{lemma}{\textbf{Lemma}}
\newtheorem{proposition}{\textbf{Proposition}}

\newtheorem{assum}{\textbf{Assumption}}
\newtheorem{remark}{\textbf{Remark}}
\usepackage{makecell}
\usepackage{hyperref}
\usepackage{algorithm}
\usepackage{algpseudocode}

\usepackage{cellspace} 
\usepackage{graphicx}
\usepackage[table]{xcolor}
\usepackage{mathtools}

\begin{document}
\title{Feasibility Analysis and Constraint Selection \\ in Optimization-Based Controllers}

\author{ Panagiotis~Rousseas$^{1}$,~\IEEEmembership{Member,~IEEE,}  Haejoon~Lee$^{2}$,~\IEEEmembership{Graduate Student Member,~IEEE,} \\ Dimos V. Dimarogonas$^{1}$,~\IEEEmembership{Fellow,~IEEE} and Dimitra Panagou$^{2,3}$,~\IEEEmembership{Senior Member,~IEEE} 
\thanks{Haejoon Lee and Dimitra Panagou would like to acknowledge the support by the National Science Foundation (NSF) under Award Number 1942907 and the Air Force Office of Scientific Research (AFOSR) under Award No. FA9550-23-1-0163.}% <-this % stops a space
% \thanks{$^\dagger$Both authors have equal contribution.}
\thanks{$^{1}$Division of Decision and Control Systems,
    School of Electrical Engineering and Computer Science,
        KTH Royal Institute of Technology, Stockholm, Sweden
        {\tt\small dimos@kth.se, rousseas@kth.se}}%
\thanks{$^{2}$Department of Robotics,
        University of Michigan, Ann Arbor, MI, USA
        {\tt\small haejoonl@umich.edu}}%
\thanks{$^{3}$Department of Aerospace Engineering,
        University of Michigan, Ann Arbor, MI, USA
        {\tt\small dpanagou@umich.edu }}%
}

\maketitle

\begin{abstract}

Control synthesis under constraints is at the forefront of research on autonomous systems, in part due to its broad application from low-level control to high-level planning, where computing control inputs is typically cast as a constrained optimization problem. Assessing feasibility of the constraints and selecting among subsets of feasible constraints is a challenging yet crucial problem. In this work, we provide a novel theoretical analysis that yields necessary and sufficient conditions for feasibility assessment {of linear constraints and based on this analysis}, we develop novel methods for feasible constraint selection {in the context of control of autonomous systems}. Through a series of simulations, we demonstrate that our algorithms achieve performance comparable to state-of-the-art methods while offering improved computational efficiency. Importantly, our analysis provides a novel theoretical framework for assessing, analyzing and handling constraint infeasibility.
% Constrained control and optimization-based controllers are at the forefront of autonomous systems' research with wide applications, from safe control to spatio-temporal task planning. In this context, assessing feasibility and selecting among subsets of feasible constraints is a challenging yet crucial problem. In this work, we provide a novel theoretical analysis that yields necessary and sufficient conditions for feasibility assessment {of linear constraints and based on this analysis}, we develop novel methods for feasible constraint selection {in the context of control of autonomous systems}. Through a series of simulations, we demonstrate that our algorithms achieve performance comparable to state-of-the-art methods while offering improved computational efficiency. 
\end{abstract}

% \begin{IEEEkeywords}
% IEEE, IEEEtran, journal, \LaTeX, paper, template.
% \end{IEEEkeywords}

% \fi
\IEEEpeerreviewmaketitle

\section{Introduction}

{
Enabling recursive feasibility of optimization-based controllers, and principled relaxation of constraints, is an essential capability towards safe and trustworthy autonomous systems; infeasibility in the face of multiple constraints can significantly compromise a system's performance and/or cause catastrophic failures. Hence, ensuring the safety, robustness and effectiveness of constrained control algorithms requires the development of novel methodologies that assess the feasibility of the underlying optimization problems, and that guide the relaxation of the constraints if the optimization problem is deemed infeasible.}
\par
The problem of constrained optimization has been extensively studied in the literature \cite{book_optim}. Finding the maximal subset of feasible constraints \cite{chinneck} --termed the maxFS problem-- is known to be NP-hard; hence approximate solutions and heuristics are usually employed \cite{782859,Chinneck2010-qg}. A common method for addressing the above is through the addition of relaxation (slack) variables to the constraints, accompanied by appropriately modifying the cost function \cite{Chinneck1996}. However, this results in a significant increase in the problem's size in the presence of multiple constraints and in general does not guarantee that the maximal feasible set is found. 
% {NP-hardness is also explored in the case of prioritized constraints  \cite{Dubois1996}, where the authors employ fuzzy logic to model constraints of different importance and propose}
% When multiple constraints in the presence of priority specifications, the authors in propose a {a fuzzy modeling approach }. 
\par
In the context of autonomous systems, a variety of methodologies aim to handle constraints. For low-level control under constraints, techniques such as Model Predictive Control (MPC) \cite{Schwenzer2021}, Reference Governors \cite{GARONE2017306} and Control Barrier Functions (CBFs) \cite{ames2016control} enforce the satisfaction of constraints either over finite horizons or pointwise in time and state; however, guaranteeing recursive feasibility remains a challenge. 
{
CBFs in particular have seen extensive development for safety, control synthesis and spatio-temporal task planning \cite{8404080}, implemented through optimization-based controllers. Most commonly, such controllers concentrate on solving a Quadratic Program (QP), termed CBF-QP controllers. However, studying and guaranteeing the feasibility of the underlying controllers remains largely unaddressed. More recently, heuristics for constraint selection based on Lagrange multipliers over dynamical system trajectories induced by optimization controllers have been proposed \cite{hardik2024}.  
}
\par 
For high-level path planning under constraints, or more generally specifications (encoded e.g., via temporal logics), the case of infeasible/incompatible specifications and how to minimally relax or violate them has also been considered \cite{Rahmani2020WhatTD, 7989177, Wongpiromsarn2020MinimumViolationPF, 6579837,10342094, 10.1145/3450267.3450542}. However, such techniques {address the high-level specifications or are specific to the application at hand, while} typically neglecting constraints at the low-level, i.e., due to dynamics, sensing and actuation limitations.
{
For instance, \cite{Rahmani2020WhatTD} treat Linear Temporal Logics (LTL) with no consideration of dynamical models, while \cite{7989177, Wongpiromsarn2020MinimumViolationPF} focus on specific road vehicle applications through sampling-based methods. They employ the notion of ``level of unsafety'' which they minimize to extract solutions in face of soft constraints. The authors in \cite{6579837} employ discrete transition models and automata, while \cite{10342094, 10.1145/3450267.3450542} use constraint relaxation via additional variables. Interestingly, \cite{10342094} follows a set-theoretic approach, instead of point-wise relaxation, but nevertheless focuses on two constraints. 
}
\par 
{More relevant to this work, in \cite{Molnar,cohen2025compatibilitymultiplecontrolbarrier} the authors investigate conditions for CBF-based constraints to be collectively feasible. Their approach can be used to guide control synthesis, with a focus on box input constraints (i.e., the $m$-dimensional input is constrained to lie within an $m$-parallelepiped). {Importantly, this enables incorporating multiple CBF constraints in a single function and can guide CBF design.} In contrast, we propose a condition {for evaluating feasibility numerically, without incorporating multiple constraints into a single formula; our condition boils down to the solution of a Linear Program (LP)}, whereas the condition in \cite{Molnar} should hold over infinitely many points (focusing on a set-theoretic approach). {This method aims at describing the properties of CBFs, and providing an algorithmic process for CBF design such that multiple constraints are satisfied, whereas our approach focuses on point-wise feasibility of a given set of constraints.}  
{
A similar approach is proposed in \cite{9993001} where feasibility checking is implemented using a relaxed LP similar to \cite{chinneck}, over a grid of points on the state-space, thus resembling the set-theoretic treatement of \cite{Molnar}. Existing methods for control also employ relaxation methods for tackling infeasibility \cite{10175590}, by incorporating relaxation variables in a QP controller. In contrast to the above, our work first, focuses on point-wise feasibility of optimization-based controllers, and second, is based on a modified LP for feasibility that does not require relaxation/slacked variables.}
\par 
% Notably, planning for control under constraints usually requires fast decision-making, often in real time. Therefore, while many modern optimization programming tools do provide feasibility analysis, these can be impractical as the scale of the problem increases, while it is desirable to avoid spending computational resources on the entire large-scale optimization problem to determine feasibility of its current instantiation. 
In our recent work \cite{rousseas2025feasibilityevaluationquadraticprograms} we developed a methodology for assessing the feasibility of QPs, which naturally arise in constrained control problems such as MPC and CBF based designs. Our approach casts the feasibility determination of the initial QP into a simpler LP, which can be solved and assessed more efficiently than the original QP problem.
% The proposed approach can then be used as a method for evaluating scenarios and determining in real-time which ones will be feasible over a horizon in the future. 
Nevertheless, several problems persist. More specifically, \cite{rousseas2025feasibilityevaluationquadraticprograms} focused only on QP-based controllers; {the technical results in \cite{rousseas2025feasibilityevaluationquadraticprograms} link (in)feasibility of QPs to (un)boundedness of the Lagrange Multipliers (LMs) associated with linear constraints. Therefore, there is no clear ``cutoff'' between feasibility/infeasibility; the set of constraints becoming infeasible when the LMs tend to infinity poses challenges to evaluating how close to infeasibility a given state is. In contrast, the method presented in the sequel proposes necessary and sufficient conditions based on the \textbf{sign} of a set of variables, providing a clear ``boundary'' between feasibility/infeasibility. 
}
% Another technical issue with using LMs for feasibility, such as in \cite{hardik2024},} is that LMs capture the interplay between the optimal solution and the constraints' set. Analyzing (in)feasibility should be based only on the structure/properties of the constraints, irrespective of the optimal solution to be extracted. 
\par 
While there exist methods that do not focus on LMs \cite{chinneck,Chinneck2010-qg}, these involve constraint relaxation via auxiliary variables; the magnitudes of the former model the degree of constraint satisfaction/violation. However, relaxing the constraints fails to decouple two critical aspects: 1) how many and which constraints are enforced, versus 2) the degree of satisfaction/violation. More precisely, such algorithms may result in solutions that enforce few constraints (i.e., minimal violation of a large number of constraints), rather than enforcing as many constraints as possible.    
\par
{Given the above, the contribution of this work is twofold:} 
\begin{enumerate}
    \item First, in contrast to existing approaches {that relax the constraints to find the maximal feasible subset}, we provide a novel analysis for extracting necessary and sufficient conditions for feasibility {over sets of linear constraints} (Sec. \ref{sec:method}), and
    \item Second, we employ the aforementioned analysis in order to provide on-line constraint selection algorithms for constrained control of dynamical systems (Sec. \ref{sec:conf_selec_methods}). 
\end{enumerate}
\par 
{The rest of the paper is organized as follows:} In Section~\ref{sec:prob}, we present some necessary preliminary notions and formulate our problem. In Sec.~\ref{sec:method}, the theoretical elements of our approach are {developed}, followed by Sec.~\ref{sec:conf_selec_methods} where two algorithms for on-line feasible constraint selection are presented. Finally, we provide simulations of the proposed algorithms in Sec.~\ref{sec:sims} with comparisons to two existing methodologies. Finally, some concluding remarks and future research directions are provided in Sec. \ref{sec:concl}. 
\par 
{
\textit{\textbf{Notation}}:
Let $\mathbb{N}$ denote the set of natural numbers, i.e., $\mathbb{N} = \{0, 1, 2, \dots\}$, $\mathbb{R}_{\geq 0}, \mathbb{R}_{< 0}$ denote a set of non-negative and strictly negative real numbers respectively. Similarly, for $n\in\mathbb{N}$, let $\mathbb{R}^n_{\geq 0}, \mathbb{R}^n_{< 0}$ denote $n$-dimensional vectors with non-negative and strictly negative real elements respectively.
Let $A\in{\mathbb{R}^{m\times c}}$ and $B\in\mathbb{R}{^c}$ denote a real-valued matrix and vector respectively. Then, $A_i\in\mathbb{R}^{m}$ denotes the matrix's $i$-th column, and $B_i \in \mathbb{R}$ denotes the vector's $i$-th element. 
Given two vectors $x \in \mathbb{R}^n,y\in\mathbb{R}^m$ let $(x,y) \in \mathbb{R}^{n+m}$ denote the stacked vector $(x,y) \triangleq [x^\top,y^\top]^\top$ and $\|x\|_p,p\in\mathbb{N}$ denotes the $p$-norm. Let $I_m$ denote the $m \times m$ identity matrix. We use ${1}_m$ and ${0}_m$ to denote $m$-dimensional vectors of ones and zeros, respectively. Let $\{0,1\}{^c}, c\in\mathbb{N}$ denote the set of vectors with $c$ elements consisting of either ones or zeros. 
}

\section{Problem Formulation}
\subsection{Preliminaries}
Consider the dynamical system:
\begin{equation}\label{eq:dynamics}
    \dot{x} = f(x) + g(x)u,
\end{equation}
where $x\in\mathbb{R}^n$ denotes the state vector, $f:\mathbb{R}^n\rightarrow\mathbb{R}^n,\ g:\mathbb{R}^n\rightarrow\mathbb{R}^{n\times m}$ are locally Lipschitz continuous functions and $u\in\mathbb{R}^m$ is the control input vector. Furthermore, for some time instance $t\in\mathbb{R}_{\geq0}$ consider {$c\in\mathbb{N}\setminus\{0\}$} affine constraints w.r.t. the input $u\in\mathbb{R}^m$ given by:
\begin{equation}\label{eq:affine_cnstr}
    A_i^\top(x,t)u \leq B_i(x,t), \quad i\in \mathcal C\triangleq \lbrace 1,\dots,c\rbrace,
\end{equation}
where $A_i : \mathbb{R}^n \times \mathbb{R}_{\geq 0}\rightarrow \mathbb{R}^m$ and $B_i: \mathbb{R}^n \times \mathbb{R}_{\geq 0}\rightarrow\mathbb{R}$ are Lipschitz continuous wrt to $x$ and $t$. %$A:\mathbb{R}^n\times \mathbb{R}_{\geq0}\rightarrow{{\mathbb{R}^{m\times c}}},\ B:\mathbb{R}^n\times \mathbb{R}_{\geq0}\rightarrow{\mathbb{R}^{c}}$ are locally Lipschitz continuous with respect to $x$ and $t$.
  %The constraints encoded {in $A,B$} consist of two subsets, namely hard constraints, which should always be satisfied, and soft constraints, which can be disregarded in case of infeasibility. 
  Without loss of generality, we assume that $n_h$ of the constraints are hard and $n_s$ of the constraints are soft, so that $n_h+n_s=c$, $\mathcal{C}_h=\{1,\dots,n_h\}$, $\mathcal C_s=\{n_h+1,\dots,c\}$. We assume:
\begin{assum}
\label{assum:hard_const_compatible}
    The set of hard constraints is always feasible, i.e., ${\forall }x\in\mathbb{R}^n, \forall t\geq 0$: 
    \begin{equation}
    \mathcal{U}_h \triangleq \left\{ u\in\mathbb{R}^m | \bigcap_{i\in\mathcal{C}_h} \{A_i^\top(x,t)u\leq B_i(x,t)\} \right\}\neq \emptyset,
    \end{equation}
    where $\mathcal{U}_h \subseteq \mathbb{R}^m$ denotes the admissible input set w.r.t. the hard constraints. %$A_i : \mathbb{R}^n \times \mathbb{R}_{\geq 0}\rightarrow \mathbb{R}^m$ denotes the $i$-th column of $A$ and $B_i: \mathbb{R}^n \times \mathbb{R}_{\geq 0}\rightarrow\mathbb{R}$ denotes the $i$-th element of $B$, where the arguments are dropped for brevity. 
\end{assum}
{
\begin{remark}
    Linear input constraints can be included in the admissible input set $\mathcal{U}_h$ \eqref{assum:hard_const_compatible} as hard constraints.
\end{remark}
}
\subsection{Modeling Disregarded Constraints}

{
\par 
{
The soft constraints are divided into two sets, namely the set of enforced constraints and the set of disregarded constraints, indexed as $\mathcal{C}_e \subseteq \mathcal{C}_s$ and $\mathcal{C}_d \subseteq \mathcal{C}_s$ respectively, such that $\mathcal{C}_e \cap \mathcal{C}_d = \emptyset$ and $\mathcal{C}_e \cup \mathcal{C}_d = \mathcal{C}_s$. The disregarded constraints are defined as those {that} are dropped from the constraint set \eqref{eq:affine_cnstr} to yield a feasible subset of constraints. In other words, {out of the soft constraints,} only the enforced ones are considered:
\begin{equation}
\notag 
    \mathcal{U}_e \triangleq \left\{ u \in \mathbb{R}^m | \bigcap_{i \in \mathcal{C}_e} \{A_i^\top(x,t)u \leq B_i(x,t)\} \right\} .
\end{equation}
Hence the admissible input set $\mathcal{U}'\subset \mathbb{R}^m$ is $\mathcal{U}'\triangleq \mathcal{U}_h \cap \mathcal{U}_e$, {whose state and time-dependence are omitted in the notation for brevity.}
}

Let the $c$-dimensional \textbf{configuration} vector $P\in \mathcal{P}$, where $\mathcal P \triangleq \{0,1\}{^c}$, be such that the $i$-th element is $0$ when the $i$-th constraint is disregarded, and $1$ when the $i$-th constraint is enforced. Denote also 
\begin{subequations}
\begin{align}
\label{eq:AB}
A(x,t) &= \bigl [A_1(x,t) | \dots | A_c(x,t)\bigr] \in \mathbb R^{m\times c},\\ B(x,t)&=\begin{bmatrix} B_1(x,t) &\dots & B_c(x,t) \end{bmatrix}^\top \in \mathbb R^c,
\end{align}
\end{subequations}
where $A_i(x,t)$ and $B_i(x,t)$, $i\in \mathcal C$ are as in  \eqref{eq:affine_cnstr}.  
\begin{remark}

The elements of the configuration vector $P\in\mathcal{P}$ corresponding to hard constraints are always set equal to one, i.e., $P_i = 1, \forall i \in \mathcal{C}_h$.
\end{remark}
Then, the set of admissible inputs $\mathcal{U}'$ can be written as a family of sets:
\begin{equation}\label{eq:input_set:config}
    \mathcal{U}'(P) = 
    \left\{
        u \in \mathbb{R}^m 
        | 
        \textrm{diag}(P)
        A^\top(x,t)u \leq
        \textrm{diag}(P)
        B(x,t)
    \right\},
\end{equation}
where {$\textrm{diag}(P)\in \{0,1\}^{c\times c}$} denotes the matrix with $P\in\mathcal{P}$ in its diagonal and zeros elsewhere, which has the effect of nullifying the disregarded constraints indexed by $\mathcal{C}_d$ in $\mathcal{U}'$. 
For any configuration $P\in\mathcal{P}$, we define an index set $\mathcal{C} \supseteq \mathcal{C}_P \triangleq \{ i \in \mathcal{C} \;| \; P_i =1 \}$ that contains the indices of the enforced constraints in $\mathcal{U}'$. 
\par
The set of \textbf{feasible configurations} $\mathcal{P}_f$ is the set that contains all configurations $P\in\mathcal P$ whose input set is non-empty, i.e.: 
\begin{equation}\label{eq:feas_conf_set_def}
\begin{gathered}
    \mathcal{P} \supseteq \mathcal{P}_f \triangleq  
    \left\{ 
            P \in \mathcal{P} \; | \;  
            \mathcal{U}'(P)\neq \emptyset 
    \right\}.
\end{gathered}\end{equation}
{Changing the configuration vector over time, effectively alters the enforced and disregarded soft constraints' sets $\mathcal{C}_d,\mathcal{C}_e$ through the input set \eqref{eq:input_set:config}.}
\par 
{
\noindent\textbf{Problem Statement:}
\label{sec:prob}
\begin{figure}[!ht]
    \centering
    \includegraphics[trim={0.75cm 3.2cm 0.5cm 1.cm},clip,width=1\linewidth]{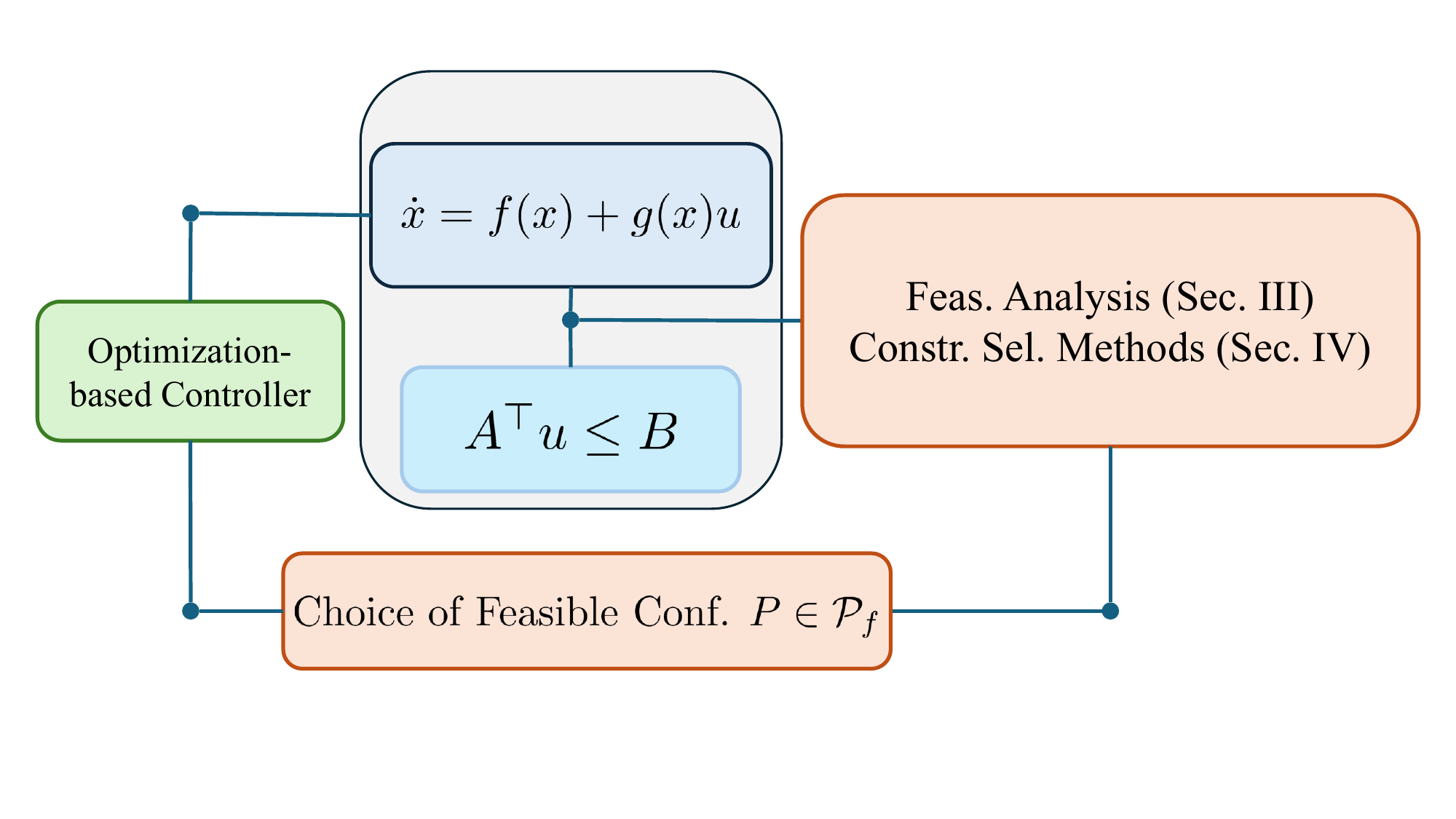}
    \caption{Overview of the Problem and its connection to the methodology (Secs. \ref{sec:method}, \ref{sec:conf_selec_methods}).}
    \label{fig:prob}
\end{figure}
Consider System \eqref{eq:dynamics} and the linear constraints \eqref{eq:affine_cnstr} as well as an optimization-based controller:
\begin{equation}\label{eq:opt:based:ctrler}
    \begin{split}
        &\quad  \qquad \qquad \dot{x} = f(x)+ g(x)u^\star(x,t), 
        \\
        &\quad  \qquad  \qquad u^\star (x,t) = \underset{u\in\mathbb{R}^m}{\arg\min}
        \left\{
            k(u)
        \right\}
        \\
        &\textrm{s.t.: } 
        \textrm{diag}(P)
        A^\top(x,t)u \leq
        \textrm{diag}(P)
        B(x,t),
    \end{split}
\end{equation}
where $k:\mathbb{R}^m\rightarrow\mathbb{R}$ is a convex function.
The goals of this work are to: 1) provide necessary and sufficient conditions for the (in)feasibility of collections of affine constraints \eqref{eq:affine_cnstr}, i.e., (non-)emptiness of the corresponding sets, and 2) 
formulate methods for on-line constraint selection for \eqref{eq:opt:based:ctrler} (see Fig. \ref{fig:prob}), i.e., choose $P(x,t)\in\mathcal{P}_f,\forall x\in\mathbb{R}^n,t\in\mathbb{R}_{\geq 0}$. 
}
}
\begin{figure*}
    \centering
    \includegraphics[trim={1.25cm 3.0cm 0.cm 3.5cm},clip,width=1\linewidth]{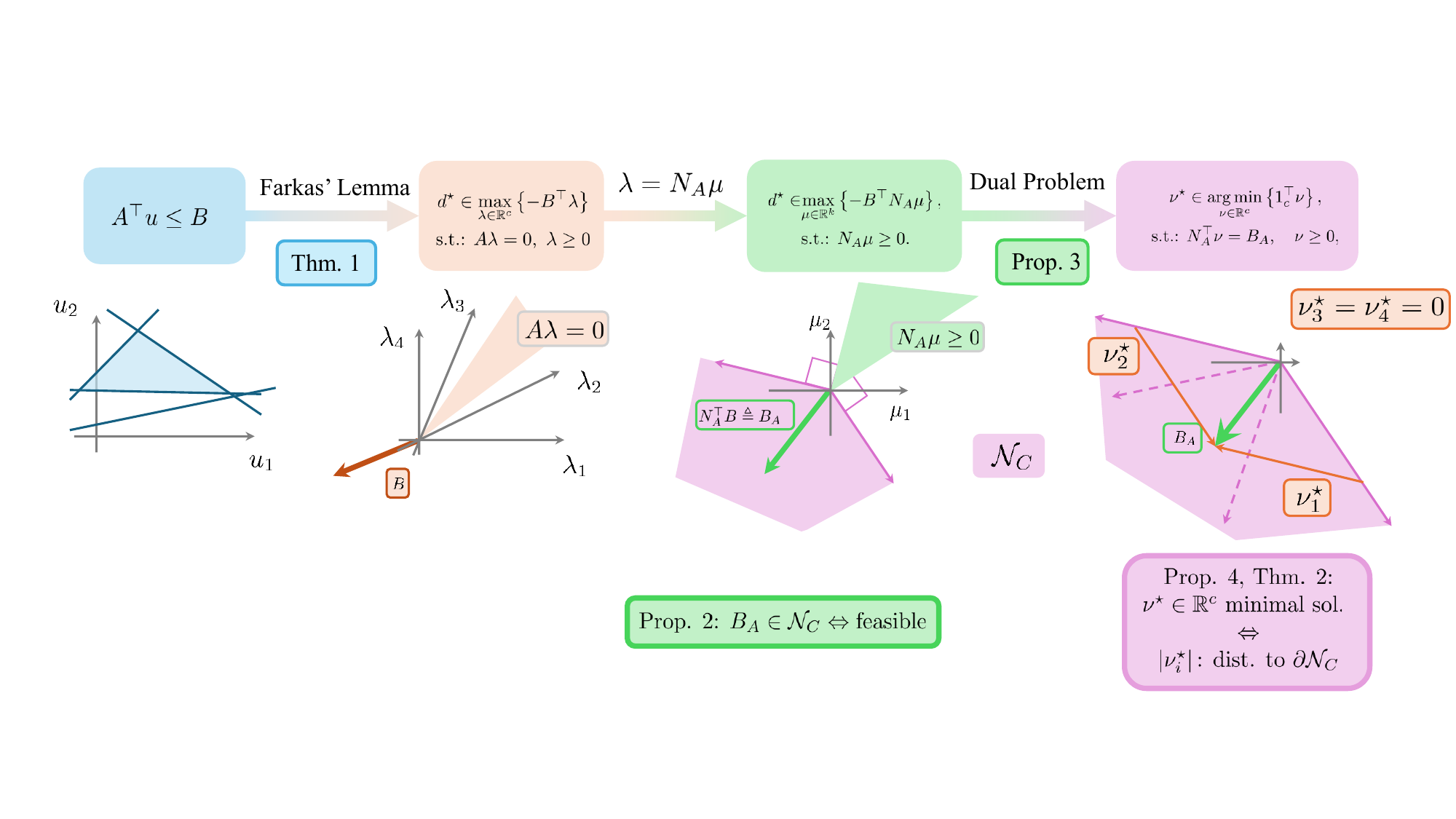}
    \caption{Overview of the Theoretical Results in Sec. \ref{sec:method}.}
    \label{fig:theory}
\end{figure*}

\section{{Feasibility Analysis}}\label{sec:method}
We now present our theoretical results for feasibility analysis {of a set of linear constraints \eqref{eq:affine_cnstr}}, outlined in Fig. \ref{fig:theory}. {This section develops the mathematical tools for feasibility analysis that are required for the online selection algorithms presented in Sec. \ref{sec:conf_selec_methods}. More specifically, we develop a criterion (denoted as $\nu^\star_{P}$ shown in \eqref{eq:polar_cone_B:star})for assessing how much each constraint contributes to the feasibility of the problem. We first analyze feasibility under the assumption that all constraints are enforced (Sec. \ref{sec:full_enforced_analysis}), and then extend the analysis to scenarios in which a subset of constraints is disregarded (Sec. \ref{sec:some_disregarded_analysis}).}

\subsection{Analysis for Enforcing all of the Constraints}\label{sec:full_enforced_analysis}

In this Subsection, we analyze the case where all of the constraints encoded in \eqref{eq:affine_cnstr} are enforced. 
% \subsubsection{\textbf{Feasibility Analysis \cite{rousseas2025feasibilityevaluationquadraticprograms}}}
{
In \cite{rousseas2025feasibilityevaluationquadraticprograms}, an LP that assesses feasibility of a set of linear constraints \eqref{eq:affine_cnstr} was proposed:
}
\begin{equation}\label{eq:feas_lp:1}
    \begin{split}
        d^\star \in \underset{\lambda\in\Lambda}{\max}
        \left\{
            -B^\top \lambda 
        \right\},        
    \end{split}
\end{equation}
where $B \triangleq  B(x,t)$, $A \triangleq  A(x,t)$ are defined in \eqref{eq:AB}, 
$\Lambda = 
    \left\{ 
        \lambda \in \textrm{null}\left(A\right) \subset \mathbb{R}{^c} 
        |
        \lambda \geq 0
\right\}$. More specifically, the set $\mathcal U=\left\{u\in\mathbb{R}^m | A^\top u \leq B \right\}$ is non-empty iff $d^\star$ in \eqref{eq:feas_lp:1} is bounded. 
In \cite{rousseas2025feasibilityevaluationquadraticprograms}, {the proof was based on the assumption that the optimization cost is of  quadratic form.} Here this result is generalized to not be limited to quadratic programs (QPs) by demonstrating its connection to Farkas' Lemma \cite{Gale1951}:
\begin{lemma}[Farkas' Lemma \cite{Gale1951}]\label{thm:farkas}
    Given two matrices $A\in{\mathbb{R}^{m\times c}},B\in\mathbb{R}{^c}$, then exactly one of the following conditions hold. (1) $A^\top u \leq B$ has a solution $u\in\mathbb{R}^m$ (Condition 1). (2) $A \lambda =0 $ (i.e., $\lambda \in  \textrm{null}(A)$) has a solution $\lambda \geq 0$ with $B^\top \lambda <0$ (Condition 2).
\end{lemma}

We begin by showing a useful proposition:
\begin{proposition}\label{prop:lambda_zero}
    If a bounded { maximum} to \eqref{eq:feas_lp:1} exists, then {the maximizer is equal to} $\lambda^\star = [0,\cdots,0]^\top \in \mathbb{R}{^c}$.
\end{proposition}
\begin{IEEEproof}
    Assume for the sake of contradiction that $\exists \bar{\lambda}\neq [0,\cdots,0]^\top$ such that $d^\star < \infty$ in \eqref{eq:feas_lp:1} is maximum. {The claim will be proven by contradiction, showing that if $\bar{\lambda}$ is the maximizer of \eqref{eq:feas_lp:1}, then $d^\star \rightarrow \infty$, thus contradicting boundedness. To achieve this, we scale $\bar{\lambda}$ by some real number $s\in\mathbb{R}$ to demonstrate that it can not be the true maximizer.}
    \par 
    Given $\bar{\lambda}$, define $d:\mathbb{R}_{\geq 0}\rightarrow \mathbb{R}$ given by {$d(s) \triangleq -s \cdot B^\top \bar{\lambda}$}. {Then, $\forall s \geq 0$, since $s A \bar{\lambda} = 0$ and $s\bar{\lambda} \geq 0$, it holds that  $s\bar{\lambda} \in\Lambda$.} We now distinguish two cases: 
    \par 
    1) If $-B^\top \bar{\lambda} > 0 $ then $\frac{\partial d}{\partial s} = -B^\top\bar{\lambda} >0 $ and hence $d=-sB^\top \bar{\lambda}$ grows unbounded as $s\rightarrow\infty$ {(by virtue of its derivative being negative)}, contradicting boundedness of $d^\star$. 
    \par 
    2) If  $-B^\top \bar{\lambda} \leq 0 $, then choosing $s^\star = 0 :  -s^\star B^\top \bar{\lambda}=0 \geq -B^\top \bar{\lambda}$. This shows that the maximizer is $s^\star \bar{\lambda} =  [0,\cdots,0]^\top $ {(rather than $\bar{\lambda}\neq [0,\cdots,0]^\top$)}, contradicting the original assumption. We have thus established that if \eqref{eq:feas_lp:1} is bounded, its {maximizer} is $\bar{\lambda} = [0,\cdots,0]^\top \in \mathbb{R}{^c}$ with the {maximum equal to $d^\star = B^\top [0,\cdots,0]^\top =  0$}.
\end{IEEEproof}

\begin{theorem}\label{thm:feas_lp_farkas}
    Given two matrices $A\in{\mathbb{R}^{m\times c}},B\in\mathbb{R}{^c}$, the set $\mathcal{U} = \left\{u\in\mathbb{R}^m | A^\top u \leq B \right\}$ is non-empty iff \eqref{eq:feas_lp:1} is bounded. 
\end{theorem}
\begin{IEEEproof} 
    (Sufficiency) Assume that \eqref{eq:feas_lp:1} is bounded. This implies that Condition 2 in Farkas' Lemma \ref{thm:farkas} cannot hold. To see this, if $\exists \bar{\lambda} \geq 0 : A\bar{\lambda} =0 $ and $B^\top \bar{\lambda} < 0 \Rightarrow-B^\top \bar{\lambda} > 0$, then following the arguments of Prop. \ref{prop:lambda_zero}, this contradicts the boundedness of \eqref{eq:feas_lp:1}. Therefore, Condition 1 of Farkas' lemma must hold, implying that $\mathcal{U}$ is non-empty.
    \par 
    % We will show this by contraposition, i.e., we need to show that if \eqref{eq:feas_lp:1} is not bounded, then the set $\mathcal{U}$ is empty. 
    (Necessity) Assume that $\mathcal{U}$ is non-empty, which implies that Condition 1 of Farkas' lemma holds. Since Condition 2 cannot then hold by Farkas' lemma, $\nexists \lambda \in \textrm{null}(A) : \lambda \geq 0$ and $-B^\top \lambda > 0$. Therefore, $\forall \lambda \in \textrm{null}(A), \lambda \geq 0$ it holds that $-B^\top \lambda \leq  0$. 
    {Therefore, the maximum of $-B^\top \lambda \leq  0$ is bounded, concluding the proof.}
    % Then, following the arguments of Prop. \ref{prop:lambda_zero},  $\lambda = [0,\cdots,0]^\top \in \mathbb{R}{^c}$ is a {maximizer} of \eqref{eq:feas_lp:1} and the maximum is $-B^\top \lambda  =  0$. 
\end{IEEEproof}
\begin{remark}\label{rem:feas_lp:contra}
    By contraposition and Thm. \ref{thm:feas_lp_farkas}, if \eqref{eq:feas_lp:1} is unbounded, then the set $\mathcal{U}$ is necessarily empty.
\end{remark}
{Thm.~\ref{thm:feas_lp_farkas} generalizes \cite[Thm.~1]{rousseas2025feasibilityevaluationquadraticprograms}, as the same result holds regardless of the structure of the cost function.} From Thm. \ref{thm:feas_lp_farkas}, the configuration $P\in\mathcal{P}$ is feasible {(and thus the set $\mathcal{U}'(P)$ \eqref{eq:input_set:config} is non-empty)} iff {$d^\star<\infty$ where}:
\begin{equation}\label{eq:feas_lp:2}
    \begin{split}
        d^\star \in \underset{\lambda\in\Lambda}{\max}
        \left\{
            -B^\top \textrm{diag}(P)\lambda 
        \right\},        
    \end{split}
\end{equation}
where 
$\Lambda = 
    \left\{ 
        \lambda \in \textrm{null}\left( \textrm{diag}(P)A\right) \subset \mathbb{R}{^c} 
        |
        \lambda \geq 0
\right\}$.
Otherwise, i.e., if {$d^\star$ in} \eqref{eq:feas_lp:2} is unbounded,  by Rem. \ref{rem:feas_lp:contra} $P$ is infeasible.
% \subsubsection{\textbf{Feasibility Analysis using the Polar Cone}}
\par 
Now, we employ Thm.~\ref{thm:feas_lp_farkas} in order to form a criterion for assessing the impact of each constraint towards rendering the set $\mathcal{U}$ empty.
Denote $\mathbb{N} \ni {k} \triangleq \textrm{dim}\left( \textrm{ker} \left( A\right) \right)$, and let {$\{e_1,\dots,e_c\}$ be a basis of $\textrm{null}(A)$, where $e_i \in \mathbb{R}^c, i \in \{1,\cdots,k\}$ are the basis vectors. Define $N_A = \bigl[ e_1 | e_2 | \cdots | e_k\bigr] \in {\mathbb{R}^{c\times k}}$, which we refer to as the nullspace basis}. Then {for any $\lambda \in \textrm{null}(A)$,} there exists a unique $\mu \in {\mathbb{R}^k}$ such that $\lambda = N_A \mu$. 
%Note that by expressing $\lambda \in \textrm{null}(A)$ in terms of the nullspace basis, i.e., $\textrm{null}(A)\ni\lambda = N_A \mu$, for $\mu \in {\mathbb{R}^k}$, 
Hence, \eqref{eq:feas_lp:1} can be written equivalently as:
\begin{equation}\label{eq:feas_lp:star}
    \begin{split}
        d^\star \in &\underset{\mu\in{\mathbb{R}^k}}{\max}
        \left\{
            -B^\top N_A \mu 
        \right\},        \\
        &\textrm{s.t.: } N_A \mu\geq 0.
    \end{split}
    \tag{\ref{eq:feas_lp:1}*}
\end{equation}
By Thm. \ref{thm:feas_lp_farkas}, if {$d^\star$ in} \eqref{eq:feas_lp:star} is bounded (resp. unbounded), $\mathcal{U}$ is nonempty (resp. empty). The LP \eqref{eq:feas_lp:star} is a conic linear problem. Given the cone $\mathcal{N} \triangleq \{ \mu \in {\mathbb{R}^k} | N_A \mu \geq 0\}$, we consider the \textit{polar cone}, i.e.:
\begin{equation}\label{eq:polar_cone}
    \mathcal{N}_C = 
    \left\{ 
        \nu \in {\mathbb{R}^k} | \nu^\top\mu \leq 0, \forall \mu \in \mathcal{N}
    \right\}.
\end{equation}
We now show that for {$d^\star$ in} \eqref{eq:feas_lp:star} to be bounded, the vector $- N_A^\top B \in {\mathbb{R}^k}$ needs to lie within the polar cone $\mathcal{N}_C$. 
\begin{proposition}\label{prop:cone}
    The {maximum of the} LP \eqref{eq:feas_lp:star} is bounded iff $-B_A\in \mathcal{N}_C$, with $B_A \triangleq  N_A^\top  B$. 
\end{proposition}
\begin{IEEEproof}
    (Sufficiency) Assume that $-B_A\in \mathcal{N}_C$. By \eqref{eq:polar_cone}: 
    \begin{equation}\notag 
    \begin{split}
        -B^\top N_A \mu \leq 0, &\forall \mu \in \mathcal{N}\Leftrightarrow \\
         -B^\top N_A \mu \leq 0, &\forall \mu \in {\mathbb{R}^k}| N_A \mu \geq 0.
    \end{split}
    \end{equation}
    Therefore, $\forall \mu \in {\mathbb{R}^k}| N_A \mu \geq 0$ it follows that $-B^\top N_A \mu \leq 0 \Rightarrow \max\left\{  -B^\top N_A \mu \right\}\leq 0$, and the {maximum of the} LP \eqref{eq:feas_lp:star} is bounded. 
    \par 
    (Necessity) We want to show that if $d^*<\infty$, then $-B_A \in \mathcal N_C $. We will prove this by contraposition, i.e., that if $-B_A \notin \mathcal{N}_C$, then {$d^*=\infty$}. For $-B_A \notin \mathcal{N}_C $, $\exists \bar{\mu}\in{\mathbb{R}^k}:N_A\bar{\mu}\geq 0$ such that $-B_A^\top \bar{\mu} > 0 \Leftrightarrow -B^\top  N_A \bar{\mu} > 0$.
    Now consider some $s\geq 0$ and define $d:\mathbb{R}_{\geq 0} \rightarrow \mathbb{R}$, given by $d(s) \triangleq  -s  B^\top N_A \bar{\mu} $. For $s\geq 0$ it holds that $N_A \bar{\mu} \geq 0 \Rightarrow N_A s \bar{\mu} = s N _A \bar{\mu}  \geq 0$. Therefore, the vector family $\mu = s\bar{\mu} \in {\mathbb{R}^k}$ satisfies the inequality in \eqref{eq:feas_lp:star} and $\frac{\partial d}{\partial s} =  -B^\top N_A \bar{\mu} > 0$. Therefore, $d$ grows unbounded as $s\rightarrow\infty$, implying that the {maximum of the} LP \eqref{eq:feas_lp:star} is unbounded, $d^\star=\infty$. Since  $-B_A \notin \mathcal{N}_C \Rightarrow d^\star = \infty$, this concludes the proof. 
\end{IEEEproof}

\subsubsection{\textbf{Feasibility Criterion}}
%In the previous subsection we have 
% Proposition \ref{prop:cone} establishes that for \eqref{eq:feas_lp:star} to be feasible, the vector $B_A \triangleq -N_A^\top B $ must lie in $\mathcal{N}_C$. This is equivalent to feasibility of the following LP:
{Prop.~\ref{prop:cone} establishes that $B_A = -N_A^\top B $ lying in $\mathcal{N}_C$ is a necessary and sufficient condition for~\eqref{eq:feas_lp:star} to be feasible. Now, consider the following problem:} 
\begin{equation}\label{eq:polar_cone_B}
    \begin{gathered}
    \nu^\star \in \underset{\nu\in\mathbb{R}{^c}}{\arg\min}\left\{{1}^\top_c \nu \right\}, \\
    \textrm{s.t.: } N_A^\top \nu = B_A , \quad  \nu \geq 0,
    \end{gathered}    
\end{equation}
where ${1}_c = \begin{bmatrix}1 & 1 &\cdots&1\end{bmatrix}^\top\in\mathbb{R}{^c}$. The equality $ N_A^\top \nu = B_A$ in \eqref{eq:polar_cone_B} implies that the maximizer $\nu^\star$ (if it exists) is the projection of $B_A$ {(which is the gradient of the objective function in \eqref{eq:feas_lp:star})} onto the column space of $\mathcal{N}_A$. We connect the LP \eqref{eq:polar_cone_B} to the LP \eqref{eq:feas_lp:star} through the following proposition:
\begin{proposition}\label{prop:lp_cone_feas}
    The maximum of the LP \eqref{eq:feas_lp:star} is bounded iff the LP \eqref{eq:polar_cone_B} is feasible. 
\end{proposition}
\begin{IEEEproof}
    Considering \eqref{eq:feas_lp:star}, its dual problem is:
\begin{equation}\label{eq:feas_lp:star:dual}
    \begin{gathered}
        \underset{\nu'\in\mathbb{R}{^c}}{\min}
        \left\{
            0 
        \right\},        \\
        \textrm{s.t.: } (-N_A)^\top \nu'= -B_A, \nu' \geq 0.
    \end{gathered}
\end{equation}
According to the weak duality principle, the maximizer of \eqref{eq:feas_lp:star} is bounded iff its dual \eqref{eq:feas_lp:star:dual} is feasible. Further, note that, for {$\mathbb{R}{^c} \ni \nu \geq 0 \Rightarrow {1}_c^\top\nu \geq 0$} hence \eqref{eq:polar_cone_B} is bounded from below. Therefore, since the constraints in \eqref{eq:feas_lp:star:dual} and \eqref{eq:polar_cone_B} are the same, feasibility of the former implies feasibility of the later and vice-versa, concluding the proof.  
\end{IEEEproof}

{
\begin{remark}
    In \cite[Thm. 2]{Molnar}, \cite[Lemma 1]{cohen2025compatibilitymultiplecontrolbarrier}, the authors employ constraints of the form:
    \begin{equation}
        c_i(x) + d_i(x)^\top u \ge 0, \quad i \in \mathcal{I} \subset \mathbb{N},\notag 
    \end{equation}
    where $c_i : \mathbb{R}^n \to \mathbb{R}$ and $d_i : \mathbb{R}^n \to \mathbb{R}^m$. The constraints are compatible for
 (1) at
    $x \in \mathbb{R}^n$ if and only if \textbf{\textit{for all}} $\lambda_i \ge 0$ the following holds:
\begin{equation}\label{eq:molner_condition_feas}
        \sum_{i \in \mathcal{I}} \lambda_i d_i(x) = 0
        \quad \Longrightarrow \quad
        \sum_{i \in \mathcal{I}} \lambda_i c_i(x) \ge 0.
        % \notag
    \end{equation}
    Consider now a fixed state $x\in\mathbb{R}^n$. Comparing \eqref{eq:molner_condition_feas} (provided in \cite[Thm. 2]{Molnar}) to \eqref{eq:feas_lp:1}, $\sum_{i \in \mathcal{I}} \lambda_i d_i(x) = 0$ corresponds to $\lambda \in\textrm{null}(A) \Leftrightarrow \lambda = N_A \mu$ which together with the condition $\lambda_i \ge 0$ yield $\lambda \in\Lambda$, while $\lambda_i c_i(x)$ is exactly $B^\top \lambda$. Then, $\lambda_i c_i(x) \ge 0$ ensures that $\underset{\lambda \in\Lambda}{\max}\left\{-B^\top \lambda\right\}\leq 0 $ is bounded from above. Hence, the two conditions are equivalent. 
     However, note that even for a single $x\in\mathbb{R}^n$, the condition in \eqref{eq:molner_condition_feas} needs to be verified $\forall \lambda \geq 0$. In contrast, the proposed LP \eqref{eq:feas_lp:star} offers a direct way of evaluating feasibility. 
    % However, note that \eqref{eq:molner_condition_feas} should hold $\forall \lambda \geq 0$. Hence, our formulation covers \eqref{eq:molner_condition_feas}, with the main difference that the proposed LP \eqref{eq:feas_lp:star} offers a direct way of evaluating feasibility, while per Rem. \ref{rem:B_cone_components} the solution to the LP can also be used to assess whether a constraint renders the problem infeasible. 
\end{remark}
}
{
\begin{remark}\label{rem:B_cone_components}
    % Each one of the components {of the solution to \eqref{eq:polar_cone_B}} $\nu^\star \in \mathbb{R}{^c}$ corresponds each of the $c$ row-vectors of the nullspace basis {in the columns of} $N_A\in{\mathbb{R}^{c\times k}}$. 
    By virtue of the equality constraint $N_A^\top \nu^\star = B_A$, each one of the components of the solution to \eqref{eq:polar_cone_B} $\nu^\star \in \mathbb{R}{^c}$ corresponds to each of the $c$ row-vectors of the matrix $N_A\in{\mathbb{R}^{c\times k}}$.
    \par 
    Note that due to the constraints in \eqref{eq:polar_cone_B}, $\nu^\star \geq 0$ element-wise. In contrast, if no such positive solution exists, (i.e., if any element of $\nu^\star$ should be negative to satisfy the equality constraint in \eqref{eq:polar_cone_B}), this implies emptiness of $\mathcal{U}$, and through Prop. \ref{prop:cone} it further implies that $B_A$ lies \textbf{outside} the polar cone \eqref{eq:polar_cone}. In case some elements of $\nu^\star $ are equal to zero, this implies that $B_A$ may lie exactly on the boundary of the polar cone and more specifically on the plane(s) formed by $(N_A)_i^\top \nu = 0, \forall i\in \{1,\cdots,c\}: \nu^\star_i = 0$.
    \par 
    Therefore, we conclude that smaller magnitudes of the elements of $\nu^\star$, $|\nu^\star_i|, i\in \{1,\cdots,c\}$ may indicate that $B_A$ is closer to the boundary of the polar cone, i.e., the distance $\underset{z \in \partial \mathcal{N}_C}{\min}\{\| z - B_A \|_2 \}$ decreases, and the set $\mathcal{U}$ tends to become empty. This is demonstrated in the last subfigure of Fig. \ref{fig:theory}, where the rows of $N_A$ are depicted with pink solid and dashed arrows and the solution for $\nu^\star$ is depicted through the orange lines, reconstructing the vector $B_A$ (green arrow). 
    % The magnitude of each of the components of $\nu^\star$ then indicates how far the vector $B_A$ is from leaving the polar cone $\mathcal{N}_A$ and hence inducing emptiness of $\mathcal{U}$. 
\end{remark}
}
To formalize the above, we prove the following: 
{
\begin{proposition}\label{prop:min-simplicial}
Let $\mathcal N = \{\mu\in\mathbb R^k \mid N_A \mu \ge 0\}$ and let its polar cone be $\mathcal N_C \coloneqq \{\nu\in\mathbb R^k \mid \nu^\top \mu \le 0,\ \forall \mu\in\mathcal N\}$. Then $\mathcal N_C$ is a pointed\footnote{That is, its rays do not extend towards infinity in all directions.}, full-dimensional polyhedral cone generated by finitely many
extreme rays $g_1,\dots,g_{\bar{c}}\in\mathbb R^k$, so that
\begin{equation}\notag 
     \mathcal N_C = \operatorname{cone}(G),\qquad
  G \coloneqq \left[ g_1 |  \dots\ | g_{\bar{c}} \right]\in\mathbb R^{k\times \bar{c}},
\end{equation}
where $\operatorname{cone}(G) \triangleq \left\{ G\lambda | \lambda \in \mathbb{R}^{\bar{c}}_{\geq 0} \right\}$.
Let $B_A\in\operatorname{int}(\mathcal N_C)$, where the latter denotes the interior of the set. Then there exists an index set
$I\subset\{1,\dots,\bar{c}\}$ with $|I|=k$ such that:
\begin{enumerate}
  \item $\{g_i\}_{i\in I}$ is linearly independent and $C_I \coloneqq \operatorname{cone}\{g_i \mid i\in I\}$
  is a $k$-dimensional simplicial cone;
  \item $B_A\in\operatorname{int} C_I$, hence there is a unique
  $\lambda_I\in\mathbb R^k_{>0}$ with $B_A = G_I \lambda_I,\ G_I \coloneqq [g_i]_{i\in I}$;
  \item the Euclidean distance from $B_A$ to the boundary of $\mathcal N_C$ coincides with the
  distance to the boundary of $C_I$:
  \begin{equation}\label{prop:cone:simplex:toprove}
      \operatorname{dist}(B_A,\partial\mathcal N_C)
    = \operatorname{dist}(B_A,\partial C_I),
  \end{equation}
\end{enumerate}
where $\operatorname{dist}(\cdot, \cdot)$ is the Euclidean distance.
Any such pair $(I,\lambda_I)$ is called a \emph{minimal simplicial representation} of $B_A$
with respect to the generator system $\{g_i\}$.
\end{proposition}
\begin{IEEEproof}
1) Since $\mathcal N_C$ is polyhedral, pointed, and full-dimensional, it has finitely many
extreme rays whose conic hull equals $\mathcal N_C$ \cite{HaynsworthFiedlerPtak1976ExtremeOperators}.
There exists a simplicial triangulation of $\mathcal N_C$ whose maximal cells are
simplicial cones $C_{I^1},\dots,C_{I^L},\ L\in\mathbb{N}$ and $I^l\subset \{ 1,\dots,\bar{c}\}, \forall l \in \{1,\dots,L\}$, each generated by $k$ linearly independent
extreme rays, and whose union equals $\mathcal N_C$, \cite{HaynsworthFiedlerPtak1976ExtremeOperators,DeLoeraHemmeckeTauzerYoshida2004Latte}.
\par
2) Because $B_A \in\operatorname{int}(\mathcal{N}_C)$, it lies in the relative interior of at least one
simplicial cone $C_{I^\ell}$ in the triangulation. Fix such an index $\ell$ and set $I\triangleq I^\ell$.
By simpliciality, the representation $B_A = G_I \lambda_I$ with unique $\lambda_I\in\mathbb R^k_{>0}$ holds.
\par
3) Each $C_{I^\ell}$ is contained in $\mathcal N_C$, so $\partial C_{I^\ell} \subset
\partial\mathcal N_C \cup\{0\}$, and hence
\begin{equation}\notag
    \operatorname{dist}(x,\partial\mathcal N_C)
  \le \operatorname{dist}(x,\partial C_{I^\ell}),\quad \forall \ell \in  \{1,\dots,L\}.
\end{equation}
Conversely, the boundary $\partial\mathcal N_C$ is the union of the facets of all simplicial
cones in the triangulation, so
\begin{equation}\notag
    \operatorname{dist}(x,\partial\mathcal N_C)
  = \min_{\ell \in  \{1,\dots,L\}} \operatorname{dist}(x,\partial C_{I^\ell}).
\end{equation}
For at least one index $\ell$ this minimum is attained; taking that $\ell$ yields item \eqref{prop:cone:simplex:toprove}.
\end{IEEEproof}

\begin{theorem}\label{thm:dist-coeff-simplicial}
Let $\mathcal N = \{\mu\in\mathbb R^k \mid N_A \mu \ge 0\}$ and let its polar cone be
$\mathcal N_C = \{\nu\in\mathbb R^k \mid \nu^\top \mu \le 0,\ \forall \mu\in\mathcal N\}$.
Assume $B_A\in\operatorname{int}\mathcal N_C$ and let $(I,\lambda_I)$ be a minimal simplicial
representation of $B_A$ in the sense of Proposition~\ref{prop:min-simplicial}, i.e.,
$|I|=k$, $G_I=[g_i]_{i\in I}\in\mathbb R^{k\times k}$ is nonsingular, and
\[
B_A = G_I \lambda_I,\quad \lambda_I\in\mathbb R^k_{>0},
\]
with
\(
\operatorname{dist}(B_A,\partial \mathcal N_C)=\operatorname{dist}(B_A,\partial C_I),
\ C_I\triangleq \operatorname{cone}(G_I).
\)
Define
\(
\nu^\star \triangleq  \lambda_I,
\nu_{\min} \triangleq  \min_{j=1,\dots,k} (\nu^\star)_j > 0,
m_I \triangleq  \sigma_{\min}(G_I),
M_I \triangleq  \max_{j=1,\dots,k}\|g_{i_j}\|_2,
\)
where $\sigma_{\min}(\cdot)$ denotes the smallest singular value.
Then the following bounds hold:
\[
m_I\,\nu_{\min}\ \le\ \operatorname{dist}(B_A,\partial\mathcal N_C)\ \le\ M_I\,\nu_{\min}.
\]
Moreover, $\nu^\star$ is the unique feasible point of the reduced LP, termed simplicial LP (over the simplex where $B_A$ lies)
\begin{equation}\label{eq:polar_cone_B:simplicial}
    \begin{gathered}
        \nu^\star = \arg\min_{\nu\in\mathbb R^k}\left\{ \mathbf{1}^\top \nu\right\}
        \\
        \quad\text{s.t.}\quad
        G_I\nu = B_A,\ \nu\ge 0.
    \end{gathered}
\end{equation}
\end{theorem}
\vspace{0.2cm}
\begin{IEEEproof}
By Proposition~\ref{prop:min-simplicial}: $\operatorname{dist}(B_A,\partial\mathcal N_C)=\operatorname{dist}(B_A,\partial C_I)$.
\\
\emph{Upper bound.}
Pick $j^\star\in\arg\min_j(\nu^\star)_j$ so that $(\nu^\star)_{j^\star}=\nu_{\min}$, and set
$\widehat{\nu}\triangleq \nu^\star-\nu_{\min}\mathbf{1}_k \ge 0$. Then the $j^\star$-th component of $\widehat{\nu}$ is zero and hence $\widehat{x}\triangleq G_I\widehat{\nu}\in\partial C_I$. Therefore,
\begin{equation}\notag
    \begin{split}
        \operatorname{dist}(B_A,\partial C_I)
        &\le \|B_A-\widehat{x}\|_2
        = \|G_I(\nu^\star-\widehat{\nu})\|_2
        \\
        &= \nu_{\min}\,\|g_{i_{j^\star}}\|_2
        \le M_I\,\nu_{\min}.
    \end{split}
\end{equation}
\\
\emph{Lower bound.}
Let $z\in\partial C_I$. Then $z=G_I\nu$ for some $\nu\ge 0$ with at least one zero entry,
so there exists $j$ with $\nu_j=0$ and hence $(\nu^\star-\nu)_j=(\nu^\star)_j\ge \nu_{\min}$,
which implies $\|\nu^\star-\nu\|_2\ge \nu_{\min}$. Since $G_I$ is non-singular $\|G_I z\|_2 \ge \sigma_{\min}(G_I)\|z\|_2,\ \forall z\in\mathbb R^k$
and therefore
\begin{equation}\notag
    \begin{split}
        \|B_A-z\|_2
        &=\|G_I(\nu^\star-\nu)\|_2
         \ge \sigma_{\min}(G_I)\,\|\nu^\star-\nu\|_2
         \\
        &\ge m_I\,\nu_{\min}.
    \end{split}
\end{equation}
Taking the infimum over $z\in\partial C_I$ yields
$\operatorname{dist}(B_A,\partial C_I)\ge m_I\nu_{\min}$, and combining with the distance
equality proves the claim.
\end{IEEEproof}

}

{
\begin{remark}
    {
        Comparing the LPs \eqref{eq:polar_cone_B} and \eqref{eq:polar_cone_B:simplicial} in Thm. \ref{thm:dist-coeff-simplicial} shows that \eqref{eq:polar_cone_B} does not necessarily guarantee that $\nu^\star$ is the minimal simplicial solution, since it includes non-extreme rays.
    }
    Nevertheless, Thm. \ref{thm:dist-coeff-simplicial} can be employed to inform which constraints (corresponding to rows of the matrix $N_A$) render $\mathcal{U}$ empty. Since $B_A \in \mathcal{N}_C$ is a necessary and sufficient condition for non-emptiness of $\mathcal{U}$, the smaller-valued elements of $\nu^\star$ \textbf{may} indicate that small variations of the corresponding constraints may cause $B_A$ to exit $\mathcal{N}_C$ {(strictly depending on whether the associated constraint forms part of the polar cone boundary -- see Thm. \ref{thm:dist-coeff-simplicial})} thus rendering $\mathcal{U}$ empty. This will be employed in Subsec. \ref{sec:meth:subsec:local_search} for local constraint selection. 
    \par 
    An example of such a minimal {simplicial} solution is depicted in the last subfigure of Fig. \ref{fig:theory}, where the non-redundant components $\nu_1^\star,\nu_2^\star > 0$ correspond to the vectors (rows of $N_A$) that form the boundary of the polar cone (shaded in pink) and $\nu_3^\star = \nu_4^\star = 0$ are the redundant components. In this example it is easy to see that solutions with $\nu_3^\star, \nu_4^\star > 0$ exist, and in such cases the corresponding magnitudes do not yield the distance of $B_A$ from the polar cone boundary. Finding the true minimal representation is a combinatorial problem, which we intend to explore in the future, based on linear (in)dependence of the rows of $N_A$.
    \par
    {
    For instance, in order to amend the issue with non-simplicial solutions, we envision two steps: 1) computing the extreme rays of the polar cone $\mathcal{N}_C$ numerically (which is however NP-hard), and 2) adding a nonlinear criterion, e.g., $\tfrac{\epsilon}{2}\| \nu\|^2$ for small $\epsilon > 0$ to the linear cost function of \eqref{eq:polar_cone_B}, resulting in the new minimization problem admitting a single solution (this would however be a QP). For $\epsilon\rightarrow0$ the QP's solution approaches the LP solution \eqref{eq:polar_cone_B}. These steps yield a 1-1 correspondence between the minimum positive element of the (modified) LP \eqref{eq:polar_cone_B} and the distance of $B_A$ to the polar cone boundary. We leave this as part of future work. 
    }
\end{remark}
}

\subsection{Feasibility Analysis under Disregarded Constraints}\label{sec:some_disregarded_analysis}
In this section we reconsider the case where only a subset of the constraints are enforced, encoded through configurations $P\in\mathcal{P}$ and their associated index sets $\mathcal{C}_P$.
Per Remark \ref{rem:B_cone_components}, the magnitude of the polar cone components of $B_A$, encoded by $\nu^\star$ in \eqref{eq:polar_cone_B} quantify how ``far'' (see the norm defined in Remark~\ref{rem:B_cone_components}) $B_A$ is from the boundary of the polar cone. Consider now a configuration $P\in\mathcal{P}_f$. Since some of the constraints have been disregarded in $P$, the following modified LP can be used to assess how far away it is from infeasibility:
\begin{equation}\label{eq:polar_cone_B:star}
    \begin{gathered}
    \nu^\star_{P} \in \underset{\nu\in\mathbb{R}{^c}}{\arg\min}\left\{\left( P - \frac{1}{2}{1}_c \right)^\top \nu \right\}, \\
    \textrm{s.t.: } N_{A}^\top \nu = B_{A}, \ \textrm{diag}\left( P\right) \nu \geq 0,
    \end{gathered}   
    \tag{\ref{eq:polar_cone_B}*}
\end{equation}
Note that for $P = {1}_c$ the equality constraint in \eqref{eq:polar_cone_B:star} are the same as in \eqref{eq:polar_cone_B}. Therefore, the former is feasible iff the latter is feasible as well. Hence, according to Thm. \ref{thm:feas_lp_farkas} and Props. \ref{prop:cone}, \ref{prop:lp_cone_feas}, $\mathcal{U}$ is nonempty iff \eqref{eq:polar_cone_B:star} is feasible. We generalize this in the following theorem:
\begin{theorem}\label{thm:conf_eval}
    A configuration $P\in\mathcal{P}$ is feasible iff \eqref{eq:polar_cone_B:star} is feasible. 
\end{theorem}
\begin{IEEEproof}
The LP \eqref{eq:polar_cone_B:star} can be recast as:
    \begin{equation}\label{eq:polar_cone_B:star:2}
    \begin{gathered}
    \nu^\star_{P}  \in \underset{\nu_e\in\mathbb{R}^{n_e},\nu_d\in\mathbb{R}^{n_d}}{\arg\min}\left\{ \frac{1}{2}{1}_{n_e}^\top \nu_e - \frac{1}{2}{1}_{n_d}^\top \nu_d \right\}, \\
    \textrm{s.t.: } N_{A,e}^\top \nu_e +  N_{A,d}^\top \nu_d 
    = B_{A}, \  \nu_e \geq 0,
    \end{gathered}   
    \tag{\ref{eq:polar_cone_B}**}
\end{equation}
where $n_e,n_d \in \mathcal{C}$ such that $n_e + n_d = c$ denote the number of enforced and disregarded constraints, respectively, and $N_{A,e}\in\mathbb{R}^{ n_e \times k}, N_{A,d}\in\mathbb{R}^{ n_d \times k}$ denote the column vectors of the nullspace basis matrix that correspond to the enforced and the disregarded constraints, respectively.
\par 
(Sufficiency:) If \eqref{eq:polar_cone_B:star} (and thus \eqref{eq:polar_cone_B:star:2}) is feasible, then there exists a $\mathbb{R}^{n_e}\ni\nu_e \geq 0$ and $\nu_d \in \mathbb{R}^{n_d}$ such that $N_{A,e}^\top \nu_e +  N_{A,d}^\top \nu_d 
    = B_{A}$. Additionally, let $\nu_d^+ \in \mathbb{R}^{n_d^+},\nu_d^-\in \mathbb{R}^{n_d^-}$ denote vectors with the positive/negative elements of $\nu_d$, with $n_d^+,n_d^- \in\mathbb{N}$ denoting the number of positive/negative elements of $\nu_d$ respectively. Then, w.l.o.g., (by rearranging the rows of $N_A$), the equality constraints in \eqref{eq:polar_cone_B:star:2} can be expressed as:
\begin{equation}\notag%\label{eq:polar_cone_B:star:3}
\begin{gathered}
   N_A^\top \nu = B_{A} \Leftrightarrow \\  
    \begin{bmatrix}
         N_{A,e}^\top, & 
        N_{A,d,+}^\top,  &
         N_{A,d,-}^\top 
     \end{bmatrix}
     \begin{bmatrix}
        \nu_e \\ 
        \nu_d^+ \\
        \nu_d^-
     \end{bmatrix} = \\
     \begin{bmatrix}
             N_{A,e}^\top, & 
             N_{A,d,+}^\top,  &
             N_{A,d,-}^\top 
         \end{bmatrix}
          \begin{bmatrix}
             B_{e}^\top, & 
             B_{d,+}^\top,  &
             B_{d,-}^\top 
         \end{bmatrix}^\top,
\end{gathered}
\end{equation}
with $\nu_e \geq 0, \nu_d^+ \geq 0, \nu_d^- < 0$, or equivalently:
\begin{equation}\label{eq:polar_cone_B:star:4}
    \begin{gathered}
        \begin{bmatrix}
             N_{A,e}^\top, & 
             N_{A,d,+}^\top,  &
             -N_{A,d,-}^\top 
         \end{bmatrix}
         \begin{bmatrix}
            \nu_e \\ 
            \nu_d^+ \\
            -\nu_d^-
         \end{bmatrix} = \\
         \begin{bmatrix}
             N_{A,e}^\top, & 
             N_{A,d,+}^\top,  &
             -N_{A,d,-}^\top 
         \end{bmatrix}
          \begin{bmatrix}
             B_{e}^\top, & 
             B_{d,+}^\top,  &
             -B_{d,-}^\top 
         \end{bmatrix}^\top 
          \\
          \Leftrightarrow
         \overline{N}_A^\top \bar{\nu} = \overline{B}_A,
    \end{gathered}
\end{equation}
with $\nu_e \geq 0, \nu_d^+ \geq 0, -\nu_d^- > 0$ and where $N_{A,e} \in \mathbb{R}^{n_e \times k}, N_{A,d,+} \in \mathbb{R}^{n_d^+ \times k}, N_{A,d,-}  \in \mathbb{R}^{n_d^- \times k}$ are formed by the columns of ${N}_A$ that correspond to $\nu_e, \nu_d^+, \nu_d^- $ while  $B_{e}  \in \mathbb{R}^{n_e}, B_{d,+}\in \mathbb{R}^{n_d^+}, B_{d,-}\in \mathbb{R}^{n_d^-}$ are formed by the elements of $B$ corresponding to the vectors $\nu_e, \nu_d^+, \nu_d^- $. This however implies that $-\overline{B}_A$ lies within the polar cone formed by the matrix $\overline{N}_A$. Note that by negating the matrix $N_{A,d,-}^\top $ and vectors $B_{d,-}, \nu_d^-$ in \eqref{eq:polar_cone_B:star:4}, this is equivalent to imposing the complementary constraints, i.e., the set:
\begin{equation}\label{eq:set_minus}
    \left\{
        u\in\mathbb{R}^m 
        \left|
            \begin{bmatrix}
                A_{e}^\top \\
                A_{d,+}^\top \\
                -A_{d,-}^\top
            \end{bmatrix} u 
            \leq 
            \begin{bmatrix}
                B_{e} \\
                B_{d,+} \\
                -B_{d,-}
            \end{bmatrix}
        \right.
    \right\}.
\end{equation}
{

To see this, note that flipping the sign of any \textbf{row} of $N_A$, is equivalent to flipping the sign of the associated Lagrange Multiplier (see Fig. \ref{fig:theory}) which is exactly equivalent to enforcing the complementary constraint, i.e., note that $\lambda = N_A \mu $ and thus $\lambda \lessgtr 0 \Leftrightarrow  N_A \mu \lessgtr 0 $.
}
Therefore, through Prop. \ref{prop:cone}, the set \eqref{eq:set_minus} is non-empty. Going back to our original definition, the aforementioned set can  be expressed as \eqref{eq:input_set:config}, which renders configuration $P$ feasible. 
\\
(Necessity:) If $P$ is feasible, then through Prop. \ref{prop:cone}, $-B_A$ belongs to the polar cone of $N_A$, which implies that a solution to:
\begin{equation}\notag
         N_{A,e}^\top
        \nu_e  = B_{A}, \nu_e\geq 0,
\end{equation}
exists, and therefore a solution to \eqref{eq:polar_cone_B:star:2} and thus to \eqref{eq:polar_cone_B:star} also exists, concluding the proof. 
\end{IEEEproof}

\section{Methods for Configuration Selection}\label{sec:conf_selec_methods}
{In the previous section, we derived the criterion $\nu_P^*$ in~\eqref{eq:polar_cone_B:star}, which quantifies the contribution of each constraint to the infeasibility of the problem, i.e., when the admissible input set satisfies $\mathcal{U}' = \emptyset$. In this section, we leverage this criterion to develop mathematically grounded methods for feasibility evaluation (Alg.~\ref{alg:feas_eval}) and for online feasible configuration selection (Alg.~\ref{alg:iter} and~\ref{alg:local_search}).}

{
\subsection{Iterative Method for Constraint Selection}
}
First the following feasibility check algorithm is proposed in Alg. \ref{alg:feas_eval}. {This algorithm utilizes the LP~\eqref{eq:polar_cone_B:star} and Thm.~\ref{thm:conf_eval} to assess the feasibility of a problem given the configuration $P$. If the LP is feasible, the algorithm returns True with $\nu_P^*$, otherwise it returns False.}
\begin{algorithm}
\caption{Feasibility Check: {\fontfamily{qcr}\selectfont FC}$(A,B,P)$}\label{alg:feas_eval}
\begin{algorithmic}[1]
\State  Given $A\in{\mathbb{R}^{m\times c}},B\in{\mathbb{R}^{c}}$, configuration $P\in\mathcal{P}$,
\State  Compute nullspace basis $N_A\in{\mathbb{R}^{c\times k}}$, $B_A= N_A^\top B$,
\State  Solve LP:
\begin{equation}\label{eq:fc:alg}
    \begin{gathered}
    \nu \in \underset{\nu\in\mathbb{R}{^c}}{\arg\min}\left\{\left( P - \frac{1}{2}{1}_c \right)^\top \nu \right\}, \\
    \textrm{s.t.: } N_{A}^\top \nu = B_{A}, \ \textrm{diag}\left( P\right) \nu \geq 0,
    \end{gathered}   
\end{equation}
\If{Problem \eqref{eq:fc:alg} is feasible}
    \State \Return TRUE, $\nu$
\ElsIf{Problem \eqref{eq:fc:alg} is infeasible}
    \State \Return FALSE
\EndIf
\end{algorithmic}
\end{algorithm}
\par 
An important observation on the LP \eqref{eq:polar_cone_B:star:2} in the context of the proof of Thm. \ref{thm:conf_eval} is that owing to minimizing the magnitude of (the positive elements of) $\nu_e\in\mathbb{R}^{n_e}_{\geq 0}$, some of the elements of the disregarded constraints' vector $\nu_d$ may become positive (i.e., $\nu_d^+ $). This is an effect of minimizing the negated sum of the elements through the term $- \frac{1}{2}{1}_{n_d}^\top \nu_d$ in \eqref{eq:polar_cone_B:star:2}; decreasing the magnitude of the elements of $\nu_e$ (accomplished by minimizing $\frac{1}{2}{1}_{n_e}^\top \nu_e$) allows for some of the elements of $\nu_d$ to become positive. This is formalized in the following corollary:
\begin{corollary}\label{col:adding_cons}
Given a feasible configuration $P\in\mathcal{P}_f$, consider the LP \eqref{eq:polar_cone_B:star} and its equivalent \eqref{eq:polar_cone_B:star:2}. Denote its minimizer as $\nu^\star = \left[ \nu_e ^\top, (\nu_d^+)^\top, (\nu_d^-)^\top \right]^\top$, where $\nu_e$ corresponds to the enforced constraints and $\nu_d^+ \in \mathbb{R}^{n_d^+}_{\geq 0}, \nu_d^- \in \mathbb{R}^{n_d^-}_{< 0}$ correspond to the disregarded constraints, with the associated constraint matrices/vectors:
\begin{equation}\notag
   \mathcal{U} =  \left\{
        u\in\mathbb{R}^m 
        \left|
            \begin{bmatrix}
                A_{e}^\top \\
                A_{d,+}^\top \\
                A_{d,-}^\top
            \end{bmatrix} u 
            \leq 
            \begin{bmatrix}
                B_{e} \\
                B_{d,+} \\
                B_{d,-}
            \end{bmatrix}
        \right.
    \right\},
\end{equation}
where $A_{e} \in \mathbb{R}^{m \times n_e}, A_{d,+} \in \mathbb{R}^{m \times  n_d^+ }, A_{d,-}  \in \mathbb{R}^{m \times n_d^-}$ are formed by the columns of $A$ that correspond to $\nu_e, \nu_d^+, \nu_d^- $ while  $B_{e}  \in \mathbb{R}^{n_e}, B_{d,+}\in \mathbb{R}^{n_d^+}, B_{d,-}\in \mathbb{R}^{n_d^-}$ are formed by the elements of $B$ corresponding to the vectors $\nu_e, \nu_d^+, \nu_d^- $. Then, enforcing the constraints 
\begin{equation}\notag
   \mathcal{U} =  \left\{
        u\in\mathbb{R}^m 
        \left|
            \begin{bmatrix}
                A_{e}^\top \\
                A_{d,+}^\top 
            \end{bmatrix} u 
            \leq 
            \begin{bmatrix}
                B_{e} \\
                B_{d,+} 
            \end{bmatrix}
        \right.
    \right\},
\end{equation}
yiels a feasible problem. 
% Enforcing the constraints corresponding to $\nu_d^+$ yields a feasible set of constraints. 
\end{corollary}
\begin{IEEEproof}
    Following the arguments of the proof of Thm. \ref{thm:conf_eval}, Eq. \eqref{eq:polar_cone_B:star:4}, positivity of $\nu_e, \nu_d^+$ implies that the vector $-\left[  B_{e}, B_{d,+}  \right]^\top$ lies within the polar cone associated with the matrix $ \begin{bmatrix}
             N_{A,e}^\top, & 
             N_{A,d,+}^\top 
         \end{bmatrix}$. Therefore, through Prop. \ref{prop:cone} enforcing the aforementioned constraints yields a feasible set.   
\end{IEEEproof}
Cor. \ref{col:adding_cons} can hence be employed to formulate an iterative method for enforcing more constraints, which is presented in Alg.~\ref{alg:iter}. In the algorithm, it starts from a feasible configuration, solves the LP \eqref{eq:feas_lp:star}, and updates the configuration by enforcing the disregarded constraints that correspond to positive components of the LP solution in \eqref{eq:feas_lp:star}. Owing to Cor. \ref{col:adding_cons}, each of the newly enforced constraints will yield feasible constraint sets.  
\par 
However, Alg. \ref{alg:iter} does not necessarily yield the maximum number of feasible constraints; note that an initial feasible configuration defines a polytope within $\mathbb{R}^m$. Since we have demonstrated that Alg. \ref{alg:iter} adds feasible constraints, its iterations can  only produce subsets within the initial polytope. Any other sets disjoint to the initial configuration's polytope, {if they exist,} cannot be obtained through this algorithm. 
\begin{algorithm}
\caption{Iterative Constraint Addition: {\fontfamily{qcr}\selectfont ICA}$(A,B,P^{(0)})$}\label{alg:iter}
\begin{algorithmic}[1]
\State  Given $A\in{\mathbb{R}^{m\times c}},B\in{\mathbb{R}^{c}}$, initial configuration $P^{(0)}\in\mathcal{P}$,
\State Check feasibility: $\mathrm{FLAG},\nu \gets${\fontfamily{qcr}\selectfont FC}$(A,B,P^{(0)})$,
\If{FLAG == FALSE}
    \State Set $P^{(0)} \gets \overline{P}$ according to the hard constraints:
    \begin{equation}
    \overline{P} = \begin{cases}
         1, & i \in \mathcal{C}_h\\
         0, & i \in \mathcal{C}_s
    \end{cases}
    \notag 
    \end{equation}
    \State Get $\nu^{(0)}$: $\mathrm{FLAG},\nu^{(0)} \gets${\fontfamily{qcr}\selectfont FC}$(A,B,P^{(0)})$,
\EndIf
\State  Split $\nu^{(0)}$ into $\left[ \nu_{e,0}^\top, (\nu_{d,0}^{+})^\top, (\nu_{d,0}^{-})^\top\right]^\top$, set $i\gets 0$,
\While{$\nu_{d,i}^{+}$ is not empty}
    \State  Update configuration:
    \begin{equation}\notag
        P^{(i+1)}_j = 
        \begin{cases}
            1,& \nu^{(i)}_j \geq 0 \\
            0,& \nu^{(i)}_j < 0
        \end{cases}, \forall j \in \{1,\cdots,C\},
    \end{equation}
    \State  Solve LP:
    \begin{equation}\notag
        \begin{gathered}
        \nu^{(i+1)} \in \underset{\nu\in\mathbb{R}{^c}}{\arg\min}\left\{\left( P^{(i+1)} - \frac{1}{2}{1}_c \right)^\top \nu \right\}, \\
        \textrm{s.t.: } N_{A}^\top \nu = B_{A}, \ \textrm{diag}\left( P^{(i+1)}\right) \nu \geq 0,
        \end{gathered}   
    \end{equation}
    \State  Split $\nu^{(i+1)}$ into $\left[ \nu_{e,i+1}^\top, (\nu_{d,i+1}^{+})^\top, (\nu_{d,i+1}^{-})^\top\right]^\top$,
    \State  $i\gets i+1$,
\EndWhile
\State \Return $P^{(i)},\nu^{(i)}$. 
\end{algorithmic}
\end{algorithm}
\par 
Alg. \ref{alg:iter} requires initialization with a feasible configuration. Note the equality constraint of \eqref{eq:polar_cone_B:star} yields:
\begin{equation}
    \begin{gathered}
        N_A^\top \nu = B_A = N_A^\top B \Leftrightarrow
        N_A^\top \left( \nu - B \right) = 0. 
    \end{gathered}
\end{equation}
Hence, for $\nu = B$ the linear constraints of the feasibility LP \eqref{eq:polar_cone_B:star} are satisfied and it is feasible as result, rendering the corresponding configuration feasible. Enforcing the constraints for which $\nu = B \geq 0$, can serve as an initialization for Alg. \ref{alg:iter}. 

\subsection{Local Search for Constraint Selection}\label{sec:meth:subsec:local_search} 
Consider again the LP \eqref{eq:polar_cone_B:star}. {As noted in Rem.~\ref{rem:B_cone_components}, the magnitude of each component of $\nu_P^*\in\mathbb{R}^c$ in~\eqref{eq:polar_cone_B:star} quantifies} the impact of each constraint towards feasibility/infeasibility, since it corresponds to a component of the vector $B_A$ within the polar cone of the matrix $N_A$. We remind the reader that $B_A$ lying within the aforementioned polar cone is a necessary and sufficient condition for the corresponding set of constraints to be feasible (see Prop. \ref{prop:cone}). Therefore, given a configuration $P\in\mathcal{P}$, the constraints can be sorted based on their corresponding value of the solution to the LP \eqref{eq:polar_cone_B:star} denoted by $\nu^\star_P$. For enforced constraints, the smaller their corresponding (positive) element of $\nu^\star_P\in \mathbb R^c$, the ``closer'' the former {may be to} rendering the problem infeasible (see Thm. \ref{thm:dist-coeff-simplicial}, Prop. \ref{prop:min-simplicial}). For disregarded constraints, the larger the corresponding (negative) element of $\nu^\star_P$ is (i.e., the closer it is to zero), the more likely it is that the constraint is feasible (see Thm. \ref{thm:dist-coeff-simplicial}, Prop. \ref{prop:min-simplicial}). 
{
To formalize this, consider Prop. \ref{prop:min-simplicial}, Thm. \ref{thm:dist-coeff-simplicial} and a change of variables similar to Thm. \ref{thm:conf_eval}. More specifically, recasting  \eqref{eq:polar_cone_B:star} as \eqref{eq:polar_cone_B:star:2} and setting $\nu_d = -\nu_d'$ yields the original LP \eqref{eq:polar_cone_B}, and Prop. \ref{prop:min-simplicial}, Thm. \ref{thm:dist-coeff-simplicial} then apply. 
}
\par 
Finally, based on the above Alg.~\ref{alg:local_search} is proposed for local constraint selection. {In this algorithm, starting from a current configuration (regardless of its initial feasibility), the algorithm iteratively modifies the enforced or disregarded sets over a user-defined search depth $D$. Alg.~\ref{alg:iter} enforces multiple additional constraints at a time, making it more computationally efficient but at the same time more susceptible to facing a "local optima" issue. In contrast, Alg.~\ref{alg:local_search} adds or disregards one constraint with the least impact on feasibility at a time. This approach helps the algorithm escape local minima and identify larger feasible configurations, though it incurs higher computational cost that scales with the search depth $D$.}

{
\subsection{Online Feasible Constraint Selection}
In this subsection the observation in Rem. \ref{rem:B_cone_components} is employed for on-line constraint selection. For system \eqref{eq:dynamics} consider the following hybrid system extension:
\begin{equation}\label{eq:dynamics:hybrid}
    \begin{split}
        \mathcal{H}
        \begin{cases}
            \begin{bmatrix}
                \dot{x} \\
                \dot{P}
            \end{bmatrix} = 
            \begin{bmatrix}
                f(x) + g(x)u \\
                {0}_C
            \end{bmatrix}, 
            & 
            x\in\mathbb{R}^n, u \in \mathcal{U}'(P),P\in\mathcal{P}_f
            \\ \\
            \begin{bmatrix}
                {x}^+ \\
                {P}^+
            \end{bmatrix} = 
            \begin{bmatrix}
                x \\
                V
            \end{bmatrix},
            & 
            x\in\mathbb{R}^n, P\notin\mathcal{P}_f, V \in \mathcal{P}_f
        \end{cases}
    \end{split}.
\end{equation}
System \eqref{eq:dynamics:hybrid} models the system evolving continuously under \eqref{eq:dynamics} under a feasible configuration $P\in\mathcal{P}_f$ \eqref{eq:feas_conf_set_def}. When infeasibility occurs, the virtual update $V\in\mathcal{P}_f$ models a change in the applied configuration denoted by $P^+$. Then the system continues to evolve according to the continuous dynamics. Let $\mathcal{T} \triangleq \bigcup_{j=0}^{J} \left([t_j,t_{j+1}] \times\{j\} \right)$ denote the hybrid time, where $j\in\{0,1,\cdots,J\},J\in\mathbb{N}$ denotes the discrete transition index and counts the number of ``jumps'' in the evolution of \eqref{eq:dynamics:hybrid}. For convenience, let any $\mathcal{T}\ni T = (t,j)$ where $t\in[t_j,t_{j+1}],j\in\{0,1,\cdots,J\}$. Then solutions to \eqref{eq:dynamics:hybrid} can be cast as hybrid arcs $x:\mathcal{T}\rightarrow\mathbb{R}^n, P:\mathcal{T}\rightarrow\mathcal{P}, u:\mathcal{T}\rightarrow\mathbb{R}^m$. 
\par 
Consequently, a decision-making scheme for $V\in\mathcal{P}_f$ needs to be designed. When infeasibility occurs, this essentially boils down to starting from the current configuration $P(T)\notin\mathcal{P}_f$ at time $T\in\mathcal{T}$ and state $x(T)$ and finding a feasible one. Formally, we model this search method as a mapping $V:\mathcal{T}\times \mathbb{R}^n \times \mathcal{P} \rightarrow \mathcal{P}_f$. Furthermore, in order to control system \eqref{eq:dynamics} an optimization-based controller can be employed. For instance, consider a QP controller:
    \begin{equation}\notag\begin{gathered}
       \pi_{A,B}(P,x,t)  = \underset{u\in\mathbb{R}^m}{\arg\min}\left\{ \|u - u_{ref}\|_2^2 \right\}\\
        \textrm{s.t.} \quad \textrm{diag}(P)A(x,t)^\top u \leq\textrm{diag}(P)B(x,t),
        \label{eq:cbf-qp}
    \end{gathered}\end{equation}
where $u_{ref}\in\mathbb{R}^m$ denotes a reference input to system~\eqref{eq:dynamics}.
We also include the polar cone components' \eqref{eq:polar_cone_B:star} evolution as the mapping $\nu^\star_P:\mathbb{R}^n\times \mathbb{R}_{\geq 0}\rightarrow\mathbb{R}^c$:
\begin{equation}\notag
    u = \pi_{A,B}(P,x,t),
\end{equation}
\begin{equation}\notag%\label{eq:polar_cone_B:star_time}
    \begin{gathered}
    \nu^\star_{P}(x,t)=  \underset{\nu\in\mathbb{R}{^c}}{\arg\min}\left\{\left( P - \frac{1}{2}{1}_c \right)^\top \nu \right\}, \\
    \textrm{s.t.: } N_{A}^\top (x,t)\nu = B_{A}(x,t), \ \textrm{diag}\left( P\right) \nu \geq 0,
    \end{gathered}   
\end{equation}
where $N_A : \mathbb{R}^n\times \mathbb{R}_{\geq 0}\rightarrow{\mathbb{R}^{c\times k}}$ denotes a nullspace basis of the matrix $A(x,t)$ and $B_A(x,t) = N^\top_A(x,t)B(x,t)$, with an overloading of notation. 
We employ two methods to achieve this task.

{
\subsubsection{\textbf{Employing Alg. \ref{alg:iter}}}
Alg. \ref{alg:iter} can be used to compute feasible configurations on-line, which requires an initial feasible configuration as an initialization. Nevertheless, owing to Assum. \ref{assum:hard_const_compatible} the hard constraints' set can be used to this end. 

\subsubsection{\textbf{Employing Alg. \ref{alg:local_search}}}
In order to acquire a feasible configuration, one/some of the enforced constraints needs to be disregarded and we remind the reader that enforced constraints correspond to positive values in the corresponding elements of the solution to \eqref{eq:polar_cone_B:star}, denoted by $\nu_P^\star$. Per Rem. \ref{rem:B_cone_components}, upon the solution of \eqref{eq:polar_cone_B:star} the constraints that correspond to the smallest (positive) elements of $\nu_P^\star$ are most likely to cause infeasibility. All of these elements are combined into Alg. \ref{alg:local_search}.
}
}

\begin{algorithm}
\caption{Local Configuration Search: {\fontfamily{qcr}\selectfont LCS}$(A,B,P,P^{-},\nu^{-},D)$}\label{alg:local_search}
\begin{algorithmic}[1]
\State  Given $A\in{\mathbb{R}^{m\times c}},B\in{\mathbb{R}^{c}}$, current configuration $P\in\mathcal{P}$, last feasible configuration $P^{-}\in\mathcal{P}_f$, LP solution vector \eqref{eq:polar_cone_B:star} $\nu^-$, search depth $\mathbb{N}\ni D\leq c$,
\State Check feasibility: $\mathrm{FLAG},\nu \gets${\fontfamily{qcr}\selectfont FC}$(A,B,P)$,
\State Initialize hard constraints configuration:
    \begin{equation}
    \overline{P}_i = \begin{cases}
         1, & i \in \mathcal{C}_h\\
         0, & i \in \mathcal{C}_s
    \end{cases},
    \notag 
\end{equation}
\If{ FLAG  = TRUE}
    \State Enforce additional constraints: 
     \begin{equation}
         P^+,\nu^+ \gets \textrm{\fontfamily{qcr}\selectfont ICA}(A,B,P),\ \notag 
     \end{equation}
     \State Save temporary variable: $P^{t}\gets P^+$,
     \For{$i = 1:D$}
        \State Split $\nu^+$ into $\left[ \nu_{e}^\top, \nu_{d}^\top\right]^\top$,
        \State 
        \parbox[t]{\dimexpr\linewidth-\algorithmicindent-10pt}{Rank elements of $\nu_{d}$ into $\nu_{d}^r \in \mathbb{R}^{c_d}_{\leq 0}$, save indices of $\nu_{d}^r$ in $\nu_{d}$ as $\{ i_1,i_2,\cdots ,i_{c_d}\}$,}
        \vspace{0.025cm}
        \For{$j = 1:c_d$}
            \State Update $P^{t}_{i_j}\gets 1$,
            \State Check feasibility:
            \begin{equation}
                \mathrm{FLAG},\nu^t \gets\textrm{\fontfamily{qcr}\selectfont FC}(A,B,P^t),\notag
            \end{equation}
            \If{FLAG = FALSE}
                \State $P^{t}_{i_j}\gets 0$
            \EndIf
        \EndFor
        \State Update $P^{+} \gets P^t,\ \nu^+ \gets \nu^t$,
     \EndFor
\ElsIf{FLAG = FALSE}
\State Save temporary variable: $P^{t}\gets P$,
    \For{$i = 1:D$}
        \State Split $\nu^-$ into $\left[ \nu_{e}^\top, (\nu_{d}^+)^\top, (\nu_{d}^-)^\top\right]^\top$,
        \State \parbox[t]{\dimexpr\linewidth-\algorithmicindent-10pt}{Rank (minimum-first) elements of $\nu_{e}$ into $\nu_{e} \in \mathbb{R}^{c_e}_{\geq 0}$, save indices of $\nu_{e}^r$ in $\nu_{e}$ as $\{ i_1,i_2,\cdots ,i_{c_d}\}$,}
        \vspace{0.025cm}
        \For{$j = 1:c_d$}
            \State Update $P^{t}_{i_j}\gets 0$,
            \State Check feasibility:
            \begin{equation}
                \mathrm{FLAG2},\nu^t \gets\textrm{\fontfamily{qcr}\selectfont FC}(A,B,P^t),\notag
            \end{equation}
            \If{FLAG2 = TRUE}
                \State $\nu^- \gets \nu^t$,
                \State \textbf{BREAK}
            \EndIf
        \EndFor
    \EndFor
    \If{FLAG2 = FALSE}
        \State\parbox[t]{\dimexpr\linewidth-\algorithmicindent-10pt}{
        \begin{equation}
             P^+,\nu^+ \gets \textrm{\fontfamily{qcr}\selectfont ICA}(A,B,\bar{P}),\ \notag 
         \end{equation}
         where $\bar{P}$ enforces only the hard constraints,}
         \vspace{0.025cm}
     \EndIf
\EndIf
\State \Return $P^{+}$. 
\end{algorithmic}
\end{algorithm}

\section{Simulations}\label{sec:sims}
In this section, the elements of the previous sections are incorporated into a scheme for efficient feasible configuration search. Herein we consider the case where the constraints are composed of hard and soft constraints indexed by $\mathcal{C}_h,\mathcal{C}_s$ respectively.

\subsection{Simulation Setup}
 Consider a model for the motion of a robot as:
$\dot x(t) = u(t)$,
where the control inputs are bounded within the set $U = [-1,1]\times [-1,1] \subset \mathbb R^2$ m/s. The robot starts at the initial position $\bar x = x(t_0)$, where $t_0 \geq \mathbb R_{\geq 0}$, and is tasked to reach a goal position $x_g$ within $T=30$ seconds, while avoiding as many undesired zones in $\mathcal C_s=\{1,\dots, n_s\}$ as possible. Each zone $i$ is modeled as a circular disk in $\mathbb R^2$ of radius $r = 1.5$ m, with centers $y_i\in \mathbb R^2$. The zones are divided into static $\mathcal C_N$ and dynamic ones $\mathcal C_D$, where $\mathcal C_N\cap \mathcal C_D=\emptyset$ and $\mathcal C_N\cup \mathcal C_D=\mathcal C_s$. Each dynamic zone $j\in \mathcal C_D$ moves with a known constant velocity $v_j \in \mathbb R^2$, maintaining constant directions and speed. We assume that the robot knows the locations athat must always be satisfiednd velocities of the zones. The zone avoidance constraints are treated as soft constraints - they may be violated if necessary during navigation - whereas the input bounds and arrival time constraints are treated as hard constraints that must always be satisfied.

Avoiding zone $i \in \mathcal C_s$ is expressed as the superlevel set of 
     $h_i(x)= ||x-y_i||_2^2 - r^2$.
Goal-reaching is guided by a Control Lyapunov Function (CLF):
 $V(x)=\|x-x_g\|^2$. Using the CBF and CLF conditions~\cite{ames2016control}, we define the affine constraint matrices $A(x,t)u \leq B(x,t)$. We first define the constraint matrices for the static zones: $A_N(x,t) = \begin{bmatrix}
    A_i(x,t)
\end{bmatrix}_{i \in \mathcal C_N}$ and $B_N(x,t) = \begin{bmatrix}
    B_i(x,t)
\end{bmatrix}_{i \in \mathcal C_N}$, where $A^\top_i(x,t) = -\frac {\partial h_i}{\partial x}$ and $B_i(x,t) = \alpha_i(h_i(x))$ for $i\in \mathcal C_N$. Similarly, we define the constraint matrices for the dynamic zones: $A_D(x,t) = \begin{bmatrix}
    A_j(x,t)
\end{bmatrix}_{j \in \mathcal C_D}$ and $B_D(x,t) = \begin{bmatrix}
    B_j(x,t)
\end{bmatrix}_{j \in \mathcal C_D}$, where $A^\top_j(x,t) = -\frac {\partial h_j}{\partial x}$ and $B_j(x,t) = \alpha_j(h_j(x))+\frac {\partial h_j}{\partial y_j} v_j$ for $j\in \mathcal C_D$. Lastly, we define the matrices for the hard constraints: $A_H(x,t)=\begin{bmatrix}
    \frac {\partial V}{\partial x}^\top & I_m & -I_m
\end{bmatrix}$ and $B_H(x,t) = \begin{bmatrix}
     -\alpha(V(x)) & \mathbf 1_m^\top& \mathbf 1_m^\top
\end{bmatrix}$. Thus, we design a controller in the form of~\eqref{eq:cbf-qp} with
where $A(x,t) = [A_N(x,t) ,A_D(x,t), A_H(x,t) ]$ and $B(x,t) = [B_N(x,t), B_D(x,t), B_H(x,t)]^\top$.

For all simulations, the initial positions of the zones, as well as the initial and goal positions of the robot, were uniformly sampled within $[-10,10]\times [-10,10]$. To ensure meaningful navigation tasks, the initial and goal positions were sampled at least $7$ m apart and lie outside the static obstacle regions. Each dynamic zone $j \in \mathcal{C}_d$ was given a random constant velocity $v_j \in \mathbb{R}^2$ with $||v_j||=2$. The initial conditions were sampled to satisfy all (both hard and soft) constraints.

We compare our methods (Alg.~\ref{alg:iter} and~\ref{alg:local_search}) against two other baselines - Baseline~1 as Chinneck's algorithm~\cite[Algorithm~1]{chinneck} and Baseline~2 as~\cite[Algorithm~2]{hardik2024}. Note that Baseline~2 does not reintroduce the constraints once they are disregarded. However, for a fair comparison in this simulation, we modified it so that the disregarded constraints are allowed to be reintroduced by assigning their corresponding Lagrange multiplier's values from the previous step to $0$. Both baseline algorithms evaluate feasibility of a given problem using a slacked linear program (LP). We evaluate performance by comparing computation times and the percentages of disregarded soft constraints across different algorithms. The simulations are run on a computer with an Intel(R) Core(TM) i5-10500H CPU @ 2.50GHz and 8.00 GB RAM.

\subsection{Varying Number of Obstacles}

We first performed a set of simulations with a varying number of the number of zones (soft constraints) and a fixed set of $5$ hard constraints corresponding to the CLF and input bounds. The number of zones ranged from $2$ to $100$ in increments of $2$, equally divided between dynamic and static constraints, across $50$ randomly generated environments. We set the search depth $D=5$ for Alg.~\ref{alg:local_search}. 

The results are shown in Fig.~\ref{fig:varying_num_obs}. In general, both of our proposed algorithms achieve performance comparable to Baseline~1 while significantly outperforming Baseline~2 in terms of disregarding constraints. However, our algorithms are in average slightly slower than Baseline~1. Algorithm~\ref{alg:iter} is faster than Baseline~$2$, while Algorithm~\ref{alg:local_search} exhibits similar average computations times to those of Baseline~2. Importantly, both proposed algorithms consistently maintain low computation times across all problem sizes, remaining below $0.07$ seconds even in the largest scenarios. In contrast, baseline algorithms, especially Baseline~1, experience frequent spikes in computation time as shown in Fig.~\ref{fig:varying_num_obs}(b), which may limit their reliabilities for real-time deployment.

\subsection{Fixed Number of Obstacles}

To further examine the performance of the proposed algorithms, we separately conducted simulations with a fixed number of $100$ obstacles across $50$ randomly generated environments. For a more comprehensive evaluation, we varied the search depth $D$ in Alg.~\ref{alg:local_search} among $1$, $5$, and $10$, and included Baseline~1 and Baseline~2 algorithms with QP feasibility checks. The distributions of disregarded-constraint percentages at each time $t$ for the different algorithms are shown in Fig.~\ref{fig:distribution}. The results are also summarized in Table~\ref{tab:records2}, showcasing the average and maximum computation times, as well as the average and maximum constraint-disregarding rates, for all considered algorithms.

 These results indicate that Alg.~\ref{alg:iter} and Alg.~\ref{alg:local_search} with different search depths achieve lower constraint-disregarding rates and thus better performance at the cost of only a slight increase in computation time compared to the baseline algorithms with LP feasibility checks, while significantly outperforming the baselines that rely on QP feasibility checks. Furthermore, similar to the previous set of simulations, our algorithms have lower maximum computation times compared to the baseline algorithms, {which is consistent with the results shown in the previous subsection.}

In summary, the simulation results demonstrate that both of proposed algorithms (Alg.~\ref{alg:iter} and~\ref{alg:local_search}) tend to achieve comparable or even better performance against the baseline algorithms in terms of the number of disregarded constraints. Compared to Baseline~1, the proposed algorithms require a slight increase in computation time, but exhibit greater reliability by consistently maintaining a significantly lower maximum computation times across a wide range of problem sizes. Moreover, unlike the heuristic baseline approaches, our methods are grounded in Cor.~\ref{col:adding_cons} and Thm~\ref{thm:dist-coeff-simplicial}, which provide a theoretically-driven criterion quantifying the impact of each constraint on feasibility/infeasibility of a problem. Finally, the proposed formulation solves LPs with only $m$ decision variables, instead of the $n_s + m$ variables required by baseline algorithms with slack variables. This reduction enables more graceful scaling with the number of soft constraints and explains the lower maximum computation times observed in our results.

\begin{figure}
    \centering
    \includegraphics[trim={1.0cm 1.75cm 1.25cm 0.0cm},clip,width=1\linewidth]{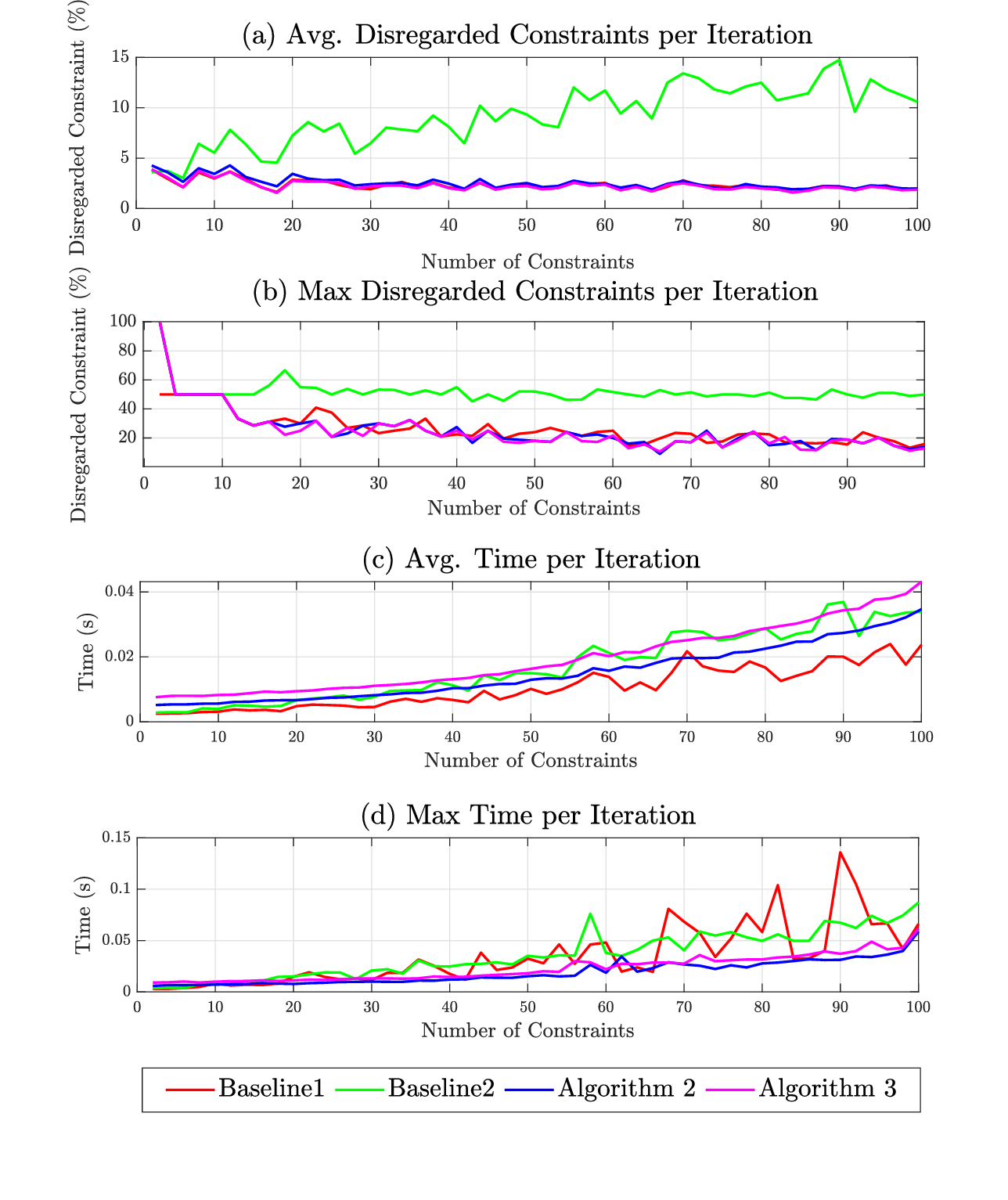}
    \caption{ Average (a) and Maximum (b) percentages of disregarded constraints, Average (c) and Maximum (d) computation times for Baseline~1, Baseline~2, Alg.~\ref{alg:iter}, and Alg.~\ref{alg:local_search} at each time step~$t$.}
    \label{fig:varying_num_obs}
\end{figure}
\begin{figure}
    \centering
    \includegraphics[trim={1.0cm 0.9cm 1.5cm 1.cm},clip,width=\linewidth]{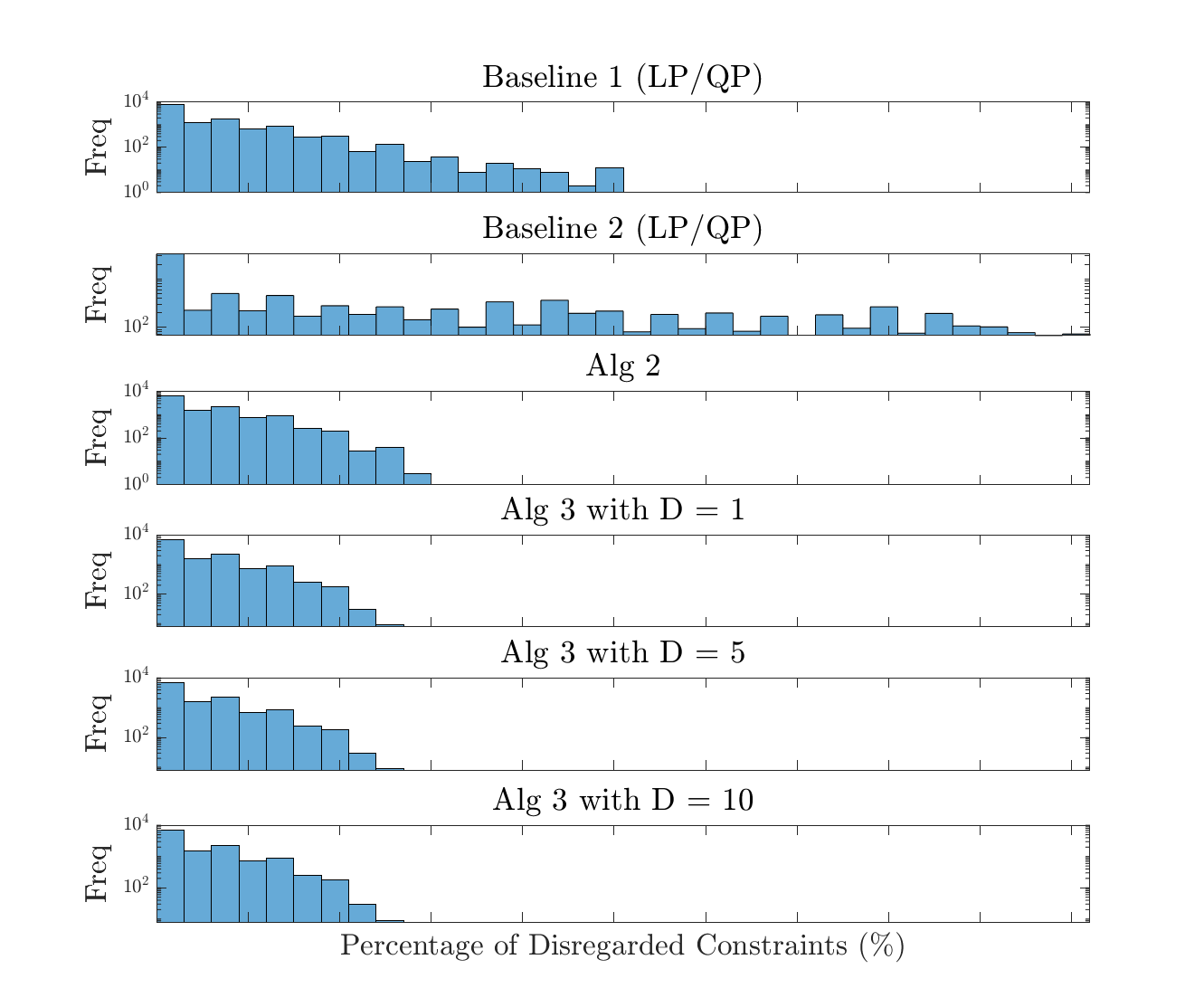}
    \caption{Histograms showing the distribution of the percentage of soft constraints disregarded (out of $50$) at each time instance across $50$ simulations for each algorithm. Frequency (in log-scale) indicates the total number of time instances at which the corresponding percentages of constraints are disregarded across $50$ different simulations.}
    \label{fig:distribution}
\end{figure}

\begin{table}[!ht]
\caption{Comparison with baseline algorithms}
\label{tab:records2}
    \centering
    \resizebox{0.95\linewidth}{!}{
    \begin{tabular}{|Sl||Sc|Sc|Sc|Sc|Sc|}
    \hline & \makecell{Avg. Comp. \\Time (s)} & \makecell{Max Comp. \\ Time (s)} &  \makecell{Avg.  \\ Drop (\%)} &  \makecell{Max \\ Drop (\%)} \\\hline\hline 
    \makecell{Baseline~$1$ \\ (LP/QP)} &  0.036 /  0.598 & 0.175 / 1.372&  2.245 &     29.000 \\\hline
    \makecell{Baseline~$2$ \\ (LP/QP)} & 0.048 / 0.286
 &   0.091 / 0.403 &  13.580 & 51.000 \\\hline\hline
    \makecell{Alg.~\ref{alg:iter}} & 0.038  &  0.050 & 2.144 & 14.000 \\\hline
   \makecell{Alg.~\ref{alg:local_search} \\ ($D=1$)} & 0.048 & 0.067
 &  2.024 & 15.000 \\\hline
 \makecell{Alg.~\ref{alg:local_search} \\ ($D=5$)} &  0.048    &   0.073 &  2.024 & 15.000 \\\hline
 \makecell{Alg.~\ref{alg:local_search} \\ ($D=10$)} & 0.049  &  0.075
 & 2.024 & 15.000 \\\hline

\end{tabular}}
\end{table}

\section{Conclusion}\label{sec:concl}
In this paper, a theoretical analysis for feasibility assessment of affine constraints was first presented and then employed to formulate on-line feasible constraint selection algorithms. The proposed framework is demonstrated to provide similar results to state-of-the art competing approaches, while exhibiting more consistent and reduced maximum computation times. 
{While the theoretical analysis was employed in this paper for online constraint selection, its importance goes beyond this problem; we believe that our necessary and sufficient conditions can be further used for tackling a number of problems which we intend to pursue in the future: the maximum feasible set problem (maxFS) \cite{chinneck}, feasible space monitoring \cite{hardik2024}, creating feasibility-aware/informed controllers and exploring adaptive techniques to render infeasible constrained control problems feasible.}

\bibliographystyle{IEEEtran}
\bibliography{main}

\begin{IEEEbiography}[{\includegraphics[width=1in,height=1.25in,clip,keepaspectratio]{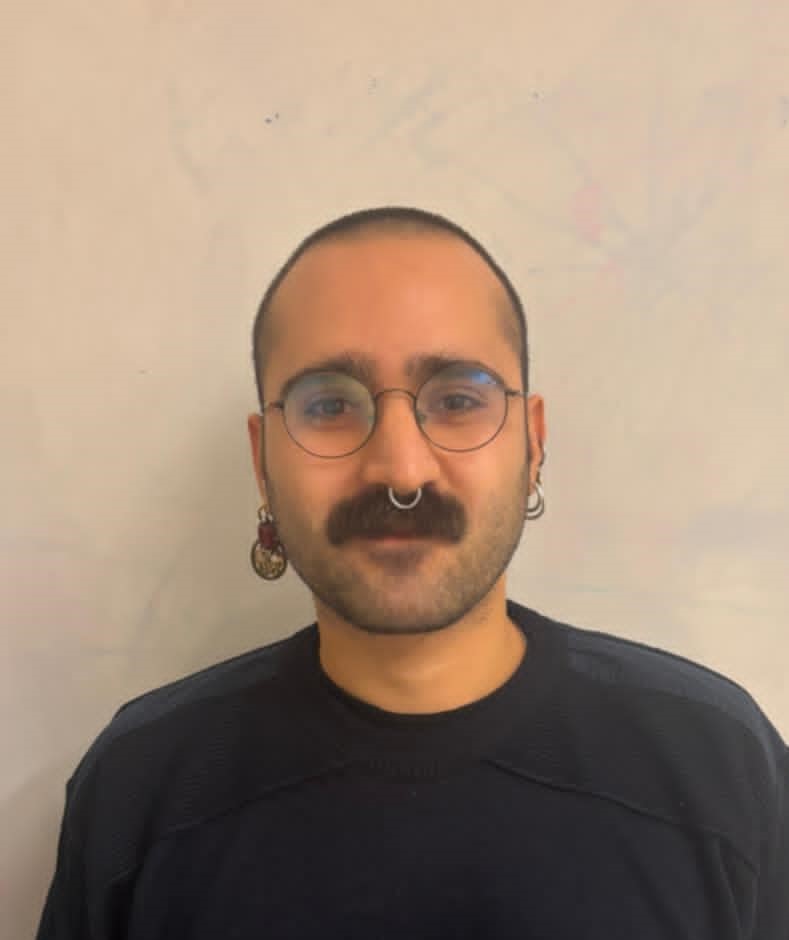}}]{Panagiotis Rousseas} (Member, IEEE) is a Post Doctoral researcher in the Royal Institute of Technology (KTH), Stockholm Sweden working in the Division of Decision and Control Systems in the group of Prof. Dimos Dimarogonas. He received his PhD in Robotics and Control at the Control Systems Laboratory of the School of Mechanical Engineering, National Technical University of Athens (NTUA), under the supervision of Professor Kostas Kyriakopoulos in 2025. He graduated with the Diploma of Mechanical Engineering from NTUA in 2020. He has visited the Robotics Department of the University of Michigan as a short term scholar in collaboration with Professor Dimitra Panagou. His PhD thesis titled "Optimal Motion Planning for mobile robots using Reinforcement Learning" concentrates on combining RL with control theory to extract solutions with provable guarantees of safety, convergence and optimality. He is currently interested in safe and robust planning and control in static and dynamic environments for safety-critical systems.
\end{IEEEbiography}

\begin{IEEEbiography}[{\includegraphics[width=1in,height=1.25in,clip,keepaspectratio]{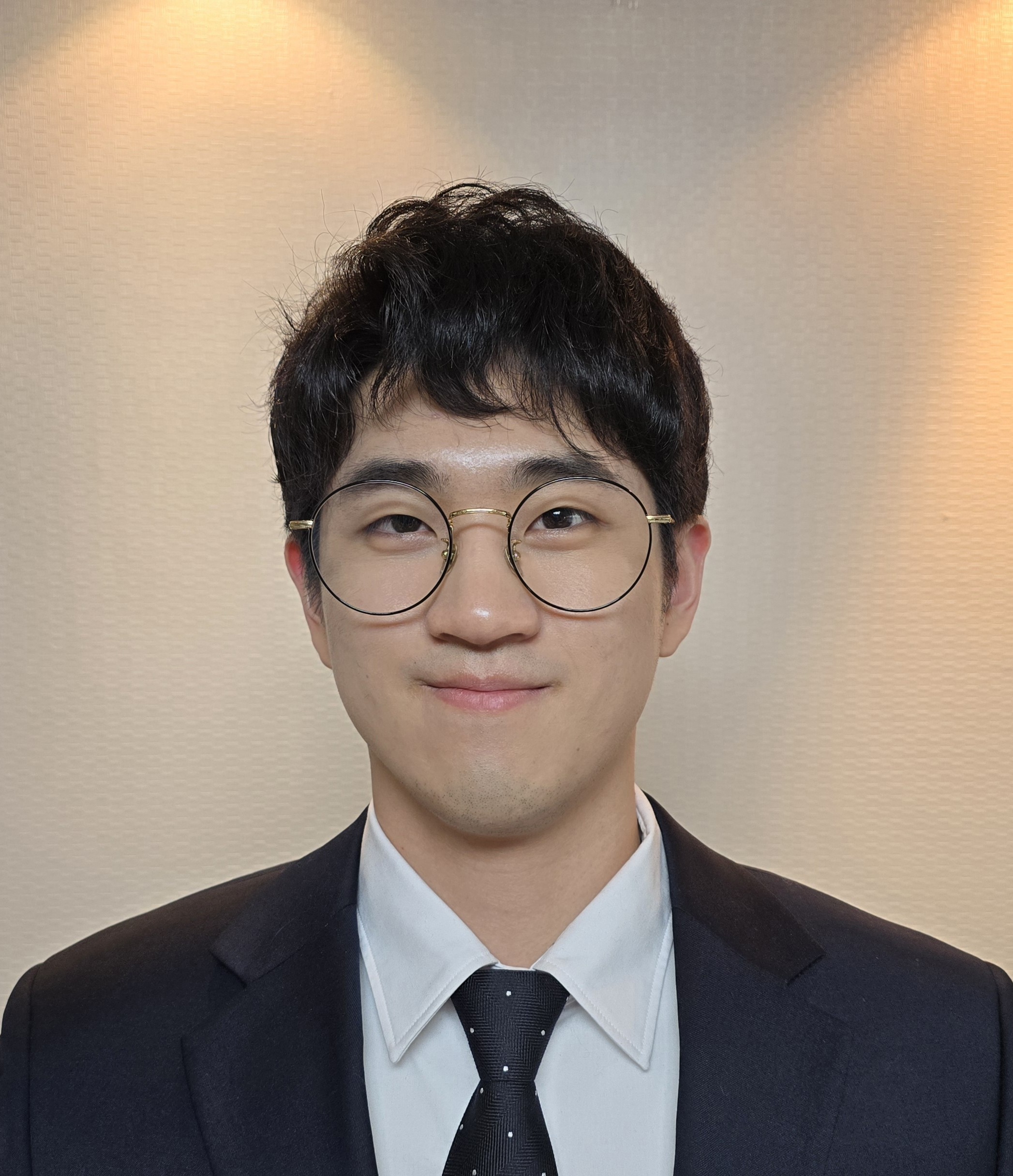}}]{Haejoon Lee} (Student Member, IEEE) received the B.S. degree in applied math and statistics from Stony Brook University, Stony Brook, NY, USA, in 2023. He earned the M.S. degree in robotics in 2025 from the University of Michigan, Ann Arbor, MI, USA, where he is currently working toward the Ph.D. degree in robotics, advised by Prof. Dimitra Panagou. His research interests include safety, resilience, and robustness of autonomous systems, with particular emphasis on distributed consensus, optimization, and learning for multi-agent systems in adversarial and uncertain environments. 
\end{IEEEbiography}

\begin{IEEEbiography}[{\includegraphics[width=1in,height=1.25in,clip,keepaspectratio]{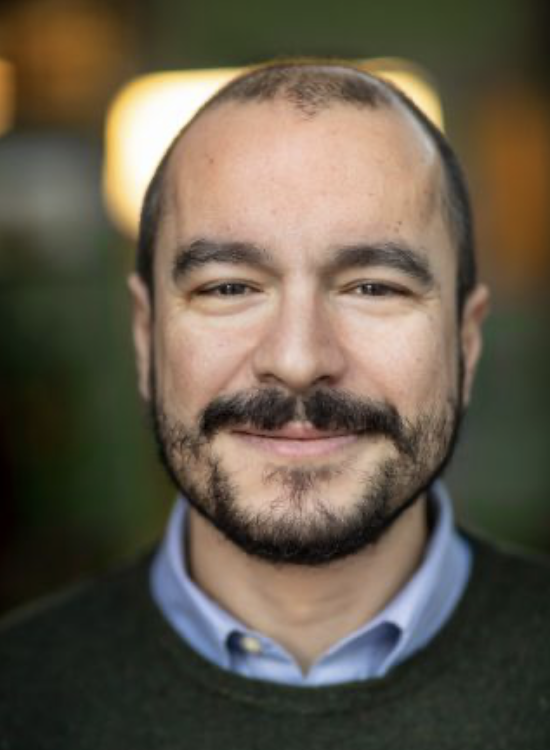}}]{Dimos V. Dimarogonas} (Fellow Member, IEEE) was born in Athens, Greece, in 1978. He received the Diploma degree in electrical and computer engineering and the Ph.D. degree in mechanical engineering from the National Technical University of Athens, Athens, Greece, in 2001 and 2007, respectively. Between 2007 and 2010, he held Postdoctoral positions with the Department of Automatic Control, KTH Royal Institute of Technology and the Laboratory for Information and Decision Systems, Massachusetts Institute of Technology, Cambridge, MA, USA. He is currently a Professor with the Division of Decision and Control Systems, School of Electrical Engineering and Computer Science, KTH Royal Institute of Technology. His current research interests include multi-agent systems, hybrid systems and control, robot navigation and manipulation, human–robot interaction, and networked control. Prof. Dimarogonas serves as an Associate Editor of Automatica and a Senior Editor of IEEE Transactions on Control of Network Systems. He was a recipient of the ERC starting Grant in 2014, the ERC Consolidator Grant in 2019, and the Knut och Alice Wallenberg Academy Fellowship in 2015. Prof. Dimarogonas is a Fellow of the IEEE.
\end{IEEEbiography}

\begin{IEEEbiography}[{\includegraphics[width=1in,height=1.25in,trim={0mm 0 0 0},clip,keepaspectratio]{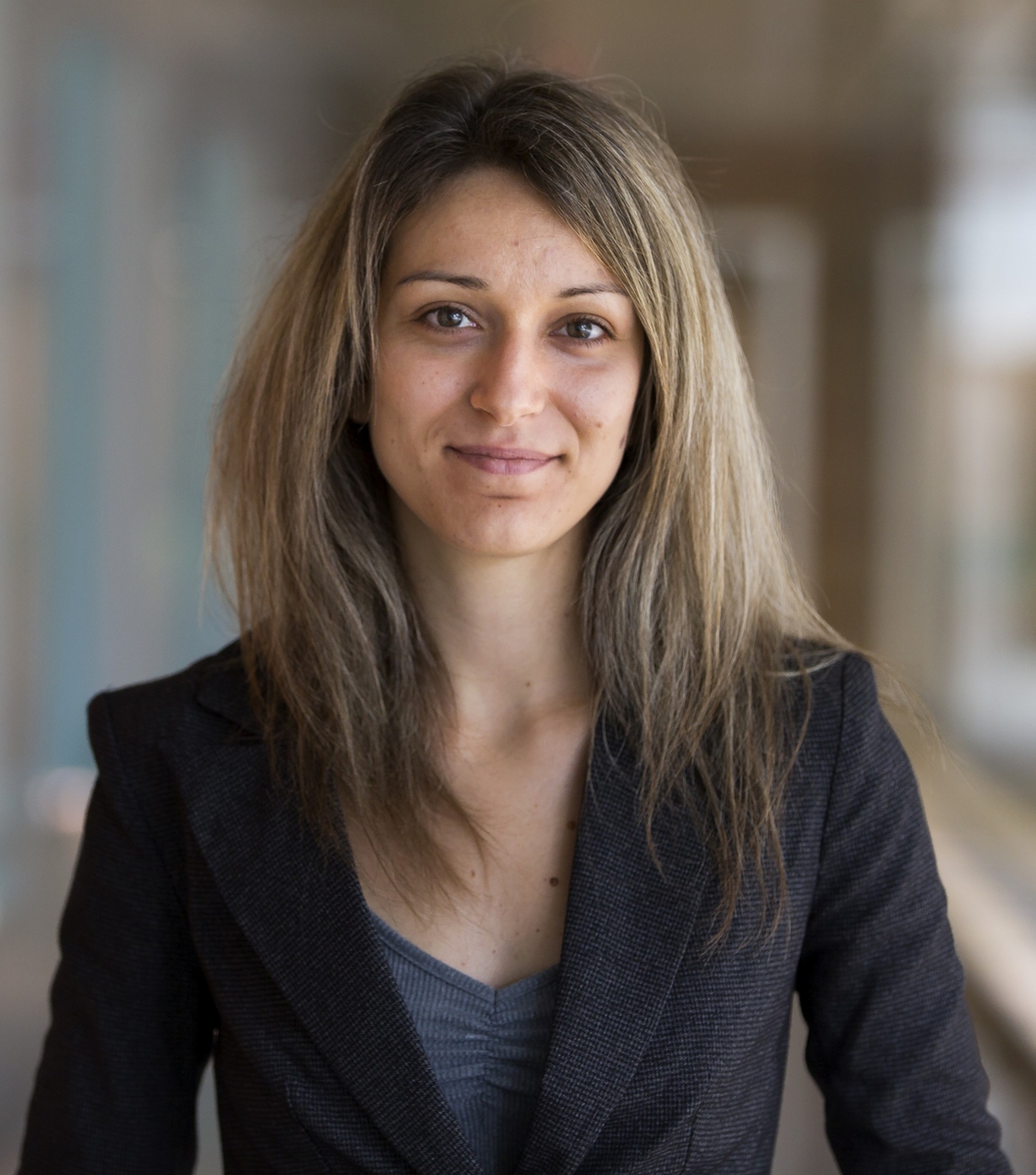}}]{Dimitra Panagou} (Diploma (2006) and PhD (2012) in Mechanical Engineering from the National Technical University of Athens, Greece) is an Associate Professor with the Department of Robotics, with a courtesy appointment with the Department of Aerospace Engineering, University of Michigan. Her research program spans the areas of nonlinear systems and control; multi-agent systems; autonomy; and  aerospace robotics. She is particularly interested in the development of provably-correct methods for the safe and secure (resilient) operation of autonomous systems with applications in robot/sensor networks and multi-vehicle systems under uncertainty. She is a recipient of the NASA Early Career Faculty Award, the AFOSR Young Investigator Award, the NSF CAREER Award, the George J. Huebner, Jr. Research Excellence Award, and a Senior Member of the IEEE and the AIAA.
\end{IEEEbiography}

\end{document}